\documentclass[10pt,a4paper,twoside,openany]{report}
\usepackage[cp1251]{inputenc} 
\usepackage{amsfonts,amsmath,amssymb,amscd}
\usepackage{bezier,graphics,graphicx}
\usepackage{epsfig,layout,latexsym}
\usepackage{theorem}
\usepackage{cite}
\usepackage{indentfirst}
\usepackage{array}

\AtBeginDocument{} \AtBeginDocument{}

\theorembodyfont{\sl}
\newtheorem{theorem}{Theorem}[section]

\newtheorem{lemma}[theorem]{Lemma}

\newtheorem{remark}[theorem]{Remark}

\numberwithin{equation}{section} \numberwithin{theorem}{section}
\newcommand{\card}{\text{card}\medskip}

\makeatletter
\renewenvironment{thebibliography}[1]
 {\section*{\centerline{\rm\textsc{Bibliography}}}%
 \@mkboth{\MakeUppercase\refname}{\MakeUppercase\refname}%
 \list{\@biblabel{\@arabic\c@enumiv}}%
 {\settowidth\labelwidth{\@biblabel{#1}}%
 \leftmargin\labelwidth
 \advance\leftmargin\labelsep
 \@openbib@code
  \usecounter{enumiv}%
  \let\p@enumiv\@empty
  \renewcommand\theenumiv{\@arabic\c@enumiv}}%
  \sloppy
  \clubpenalty4000
  \@clubpenalty \clubpenalty
  \widowpenalty4000%
  \sfcode`\.\@m
  \setlength{\itemsep}{-0.1cm}}
  {\def\@noitemerr
  {\@latex@warning{Empty 'thebibliography' environment}}%
 \endlist}
\renewcommand{\@biblabel}[1]{#1.}

\makeatother \hfuzz=0.5pt \tolerance=500

\begin{document}

\setcounter{equation}{0} \setcounter{figure}{0} \setcounter{table}{0} \setcounter{footnote}{0} \setcounter{section}{0}

\begin{center}
\textbf{About optimization of methods for mixed derivatives of bivariate functions }

\setcounter{footnote}{2} \footnotetext{\textit{Key words}. Numerical differentiation, Legendre polynomials, truncation
method, minimal radius of Galerkin information.}
\end{center}

\def\headlinetitle{Error bounds for Fourier-Legendre truncation method}

\vspace{-4mm} \noindent


\vspace*{5mm}
\centerline{\textsc { Y.V. Semenova $\!\!{}^{\dag}$, S.G. Solodky  $\!\!{}^{\dag,\ddag}$} }


%
\vspace*{5mm}
\centerline{$\!\!{}^{\dag}\!\!$ Institute of Mathematics, National Academy of Sciences of Ukraine, Kyiv}
\centerline{$\!\!{}^{\ddag}\!\!$ University of Giessen, Department of Mathematics, Giessen, Germany}

\vspace{4mm}
\begin{small}
\begin{quote}
\textsc{Abstract.}

\end{quote}
\end{small}

\begin{small}
\begin{quote}
The problem of optimal recovering high-order mixed derivatives of bivariate functions with finite
smoothness is studied. On the basis of the truncation method, an algorithm for numerical differentiation is
constructed, which is order-optimal both in the sense of accuracy and in terms of the amount of involved Galerkin
information.
\end{quote}
\end{small}

\section{Description of the problem}

Currently, many research activities on the problem of stable numerical differentiation have been taking place due to
the importance of this tool in such areas of science and technology as finance, mathematical physics, image processing,
analytical chemistry, viscous elastic mechanics, reliability analysis, pattern recognition, and many others. Among these investigations we highlight \cite{Dolgopolova&Ivanov_USSR_Comput_Math_Math_Phys_1966_Eng}, which is the first
publication on numerical differentiation in terms of the theory of ill-posed problems.
Further research \cite{Dolgopolova&Ivanov_USSR_Comput_Math_Math_Phys_1966_Eng} has been continued in numerous
publications on numerical differentiation (see, for example, \cite{Ramm_1968_No11}, \cite{VasinVV_1969_V7_N2},
\cite{Groetsch_1992_V74_N2}, \cite{Hanke&Scherzer_2001_V108_N6}, \cite{Ahn&Choi&Ramm_2006}, \cite{Lu&Naum&Per},
\cite{Nakamura&Wang&Wang_2008}, \cite{Zhao_2010}, \cite{Zhao&Meng&Zhao&You&Xie_2016}, \cite{Meng&Zhaoa&Mei&Zhou_2020},
\cite{SSS_CMAM}),  covering different classes of the functions and the types of proposed methods. Despite the abundance
of works on this topic, the problem for  recovery of high order  derivatives was considered only in a few publications,
among which we note \cite{EgorKond_1989}, \cite{RammSmir_2001}, \cite{ Wang_Hon_Ch_2006},
\cite{Qian&Fu&Xiong&Wei_2006}, \cite{Nakamura&Wang&Wang_2008} and \cite{SSS_CMAM}. In particular, the results of
\cite{SSS_CMAM} have opened  perspective for further investigation of numerical methods for recovery of high order
derivatives. Namely as the main criteria of method's efficiency  have been taken its ability to achieve  the optimal
order of accuracy by using the minimal amount of discrete information. Note that particular these aspects of numerical
differentiation  remains still insufficiently studied. The present paper continues the research of \cite{SSS_CMAM},
\cite{Sem_Sol_2021} and proposes a numerical method for recovering the high order mixed derivatives of smooth bivariate
functions. The method is not only stable to small perturbations of the input data, but is also optimal in terms of
accuracy and quantity of involved Fourier-Legendre coefficients, and also has a simple numerical implementation.

Let $\{\varphi_k(t)\}_{k=0}^\infty$  be the system of Legendre polynomials  orthonormal on $[-1,1]$ as
$$
\varphi_k(t)=\sqrt{k+1/2}(2^kk!)^{-1}\frac{d^k}{dt^k}[(t^2-1)^k], \quad k=0,1,2,\ldots .
 $$

 By $L_2=L_2(Q)$   we mean space of square-summable  on $Q=[-1,1]^2$ functions $f(t,\tau)$ with inner  product
$$
\langle f, g\rangle=\int_{-1}^{1}\int_{-1}^{1}f(t,\tau)g(t,\tau)d \tau d t
$$
and standard norm
$$
\|f\|_{L_2}^2=\sum_{k,j=0}^{\infty}|\langle f, \varphi_{k,j} \rangle|^2 < \infty ,
$$
where $$ \langle f, \varphi_{k,j}\rangle=\int_{-1}^{1}\int_{-1}^{1}f(t,\tau)\varphi_k(t)\varphi_j(\tau)d\tau dt, \quad
k,j=0,1,2,\ldots,
$$
are Fourier-Legendre coefficients of $f$. Let $\ell_p$, $1\leq p\leq\infty$, be the space of numerical sequences
$\overline{x}=\{x_{k,j}\}_{k,j\in\mathbb{N}_0}$, $\mathbb{N}_0=\{0\}\bigcup\mathbb{N}$, such that the corresponding
relation
$$
\|\overline{x}\|_{\ell_p}  := \left\{
\begin{array}{cl}
\bigg(\sum\limits_{k,j\in\mathbb{N}_0} |x_{k,j}|^p\bigg)^{\frac{1}{p}} < \infty ,
 \ & 1\leq p<\infty ,
\\\\
\sup\limits_{k,j\in\mathbb{N}_0}  |x_{k,j}| < \infty ,
  \ & p=\infty ,
\end{array}
\right.
$$
is fulfilled.

We introduce   the space of functions
$$
  L_{s,2}^\mu(Q)=\{f\in L_2(Q): \quad \|f\|_{s,\mu}^s:=\sum_{k,j=0}^{\infty} ({\underline{k\cdot j}})^{s\mu}|\langle f,
  \varphi_{k,j}\rangle|^s<\infty\},
$$
where $\mu>0$,\ $1\le s<\infty$,\ $\underline{k}=\max\{1,k\}$,\ $k=0,1,2,\dots$. Note that in the future we will use
the same notations both for space and for a unite ball from this space: $L_{s,2}^{\mu} = L_{s,2}^{\mu}(Q) = \{f\in
L_{s,2}^{\mu}\!: \|f\|_{s,\mu} \leq 1\}$, what we call  a class of functions. What exactly is meant by $
L_{s,2}^{\mu}$, space or class, will be clear depending on the context in each case. It should be noted that
$L^{\mu}_{s,2}$ is generalization of the class of bivariate functions with dominating mixed partial derivatives.
Moreover, let $C=C(Q)$ be the space of continuous on $Q$ bivariate  functions.

We represent a function $f(t,\tau)$ from $L_{s,2}^{\mu}$ as
$$
 f(t,\tau) = \sum_{k,j=0}^{\infty}\langle f, \varphi_{k,j}\rangle \varphi_k(t)\varphi_j(\tau),
$$
and by its mixed derivative we mean the following series
\begin{equation}\label{r_deriv}
f^{(r_1,r_2)}(t,\tau) =  \sum_{k,j=r}^{\infty}\langle f, \varphi_{k,j}\rangle
\varphi^{(r_1)}_k(t)\varphi^{(r_2)}_j(\tau), \quad r_1, r_2=1,2,\ldots.
\end{equation}
Assume that instead of the exact values of  the Fourier-Legendre coefficients $\langle f , \varphi_{k,j} \rangle$ only
some their perturbations are known with the error level $\delta$ in the metrics of $\ell_p$, $1\leq p\leq\infty$.
%
 More accurately, we assume that there is a sequence of numbers $\overline{f^\delta}= \{\langle
f^\delta, \varphi_{k,j} \rangle\}_{k,j\in\mathbb{N}_0}$ such that for $\overline{\xi}=
\{\xi_{k,j}\}_{k,j\in\mathbb{N}_0}$, where $\xi_{k,j}=\langle f-f_\delta,\varphi_{k,j}\rangle$, and for some $1\leq
p\leq \infty$ the relation
 \begin{equation}\label{perturbation}
   \|\overline{\xi}\|_{\ell_p} \leq \delta , \quad 0<\delta <1 ,
 \end{equation}
is true.

The research of this work is devoted to the optimization of methods for recovering the derivative (\ref{r_deriv}) of
functions from classes $L^{\mu}_{s,2}$. Further, we give a strict statement of the problem to be studied. In the
coordinate space $[r_1,\infty)\times[r_2,\infty)$ we take an arbitrary bounded domain $\Omega$. By $\card(\Omega)$ we
mean the number of points that make up $\Omega$ and by the information vector $G(\Omega, \overline{f}^{\delta})\in
\mathbb{R}^{N}$, $\card (\Omega) = N$, we take the set of perturbed values of Fourier-Legendre coefficients
$\left\{\langle f^{\delta}, \varphi_{k,j} \rangle\right\}_{(k,j)\in\Omega}$.

Let  $X=L_{2}(Q)$ or $X=C(Q)$.
 By numerical differentiation algorithm we mean any mapping $\psi^{(r_1,r_2)} = \psi^{(r_1,r_2)}(\Omega)$ that corresponds to the information vector
 $G(\Omega, \overline{f}^{\delta})$ an element $\psi^{(r_1,r_2)}(G(\Omega, \overline{f}^{\delta})) \in X$,
 which is taken as an approximation to the derivative (\ref{r_deriv}) of function $f$ from the class $L_{s,2}^{\mu}$.
We denote by  $\Psi(\Omega)$ the set of all algorithms $\psi^{(r_1,r_2)}(\Omega):\,\mathbb{R}^{N}\rightarrow X$, that
use the same information vector $G(\Omega, \overline{f}^{\delta})$.

We do not require from algorithms which belong to $\Psi(\Omega)$, in generally speaking,  either linearity or even
stability. The only condition for algorithms from $\Psi(\Omega)$ is to use input information in the form of perturbed
values of the coefficients Fourier-Legendre with indices from the domain $\Omega$  of the coordinate plane. Such the
general understanding of the algorithm is explained by the desire to compare the widest possible range of possible
methods of numerical differentiation.

The error of the algorithm $\psi^{(r_1,r_2)}$ on the class $L^{\mu}_{s,2}$ is determined by the quantity
 $$
 \varepsilon_{\delta}(L^{\mu}_{s,2}, \psi^{(r_1,r_2)}(\Omega), X, \ell_p)
 = \sup\limits_{f\in L^{\mu}_{s,2}, \, \|f\|_{s,\mu}\leq 1}
 \ \, \sup\limits_{\overline{f}^{\delta}: \, (\ref{perturbation})}
 \| f^{(r_1,r_2)} - \psi^{(r_1,r_2)}(G(\Omega, \overline{f}^{\delta})) \|_X .
 $$

The minimal radius of the Galerkin information for the problem of numerical differentiation on the class
$L^{\mu}_{s,2}$ is given by
$$
R^{(r_1,r_2)}_{N,\delta} (L^{\mu}_{s,2}, X, \ell_p) = \inf\limits_{\Omega: \, \card(\Omega)\leq N}
 \ \, \inf\limits_{\psi^{(r_1,r_2)}\in\Psi(\Omega)} \varepsilon_{\delta}(L^{\mu}_{s,2}, \psi^{(r_1,r_2)}(\Omega), X, \ell_p) .
$$

This  value  $R^{(r_1,r_2)}_{N,\delta} (L^{\mu}_{s,2}, X, \ell_p)$ describes the minimal possible accuracy in the
metric of $X$, which can be achieved by numerical differentiation of arbitrary function $f\in L_{s,2}^{\mu}$ , while
using not more than $N$ values of its Fourier-Legendre coefficients that are $\delta$-perturbed in the $\ell_p$ metric.
Note that the minimal radius of Galerkin information in the problem of recovering the first partial derivative was
studied in \cite{Sol_Stas_UMZ2022}, and for other types of ill-posed problems, similar studies were previously carried
out in \cite{PS1996}, \cite{Mileiko_Solodkii_2014}. It should be added that the minimal radius characterizes the
information complexity of the considered problem and is traditionally studied within the framework of the IBC-theory
(Information- Based Complexity Theory), the foundations of which are laid in monographs \cite{TrW} and \cite{TrWW}.

The aim of our research is to find sharp estimates (in the power scale) for the quantities $R^{(r_1,r_2)}_{N,\delta}
(L^{\mu}_{s,2}, C, \ell_p)$ and $R^{(r_1,r_2)}_{N,\delta} (L^{\mu}_{s,2}, L_{2}, \ell_p)$.

\section{Truncation method. Error estimate in $L_2$-- metric }

It should be noted that at the moment a number of approaches were developed for numerical differentiation (see, for
example, \cite{Cul71}, \cite{And84}, \cite{Qu96}, \cite{RammSmir_2001}, and also see \cite{SSS_Rew2021} and its
bibliography). All these methods are accepted to divide  into three groups (see \cite{RammSmir_2001}): difference
methods, interpolation methods and regularization methods. As is known, the first two types of methods have their
advantage in the simplicity of implementation, but they guarantee satisfactory accuracy only in the case of exactly
given input data about the differentiable function.  In the same time regularization methods give stable approximations
to the desired derivatives in the case of perturbed input data  but most of them (for example, the Tikhonov method and
its various variations) are quite complicated for numerical realization in view of their integral form and require
hard-to-implement rules for determination of regularization parameters (see \cite{RammSmir_2001}). Recently in
\cite{SSS_CMAM} a concise numerical method, called the truncation method, has been proposed as stable and simple
approach to numerical differentiation of multivariable functions. The essence of this method is to replace the Fourier
series (\ref{r_deriv}) with a finite Fourier sum using perturbed data $\langle f_\delta, \varphi_{k,j} \rangle$. In the
truncation method  to ensure the stability of the approximation and achieve the required order accuracy,  it is
necessary to choose properly the discretization parameter, which here serves as a regularization parameter. So, the
process of regularization in method under consideration consists in matching the discretization parameter with the
perturbation level $\delta$ of the input data. Simplicity of implementation is the main advantage of this method.

In the case of an arbitrary bounded domain $\Omega$ of the coordinate plane $[r_1,\infty)\times[r_2,\infty)$, the
truncation method for differentiating functions of two variables has the form
$$
\mathcal{D}_\Omega^{(r_1,r_2)} f_\delta(t,\tau)=\sum_{(k,j)\in \Omega} \langle f_\delta, \varphi_{k,j} \rangle
\varphi_k^{(r_1)}(t)\varphi_j^{(r_2)}(\tau).
$$
In order to increase the efficiency of the approach under study, we take a hyperbolic cross as the domain $\Omega$ of
the following form
$$
\Omega = \Gamma_{n,\gamma} := \{ (k,j):\ k\cdot j^\gamma \leq n,\quad k=r_1,\ldots, n/r_2^\gamma,\quad j=r_2,\ldots,
(n/r_1)^{1/\gamma}\}, \qquad \gamma \geq 1 .
$$
Then the proposed version of the truncation method can be written as
\begin{equation} \label{ModVer}
\mathcal{D}_{n,\gamma}^{(r_1,r_2)} f_\delta(t,\tau) = \sum_{k=r_1}^{n/r_2^\gamma}\, \sum_{j=r_2}^{(n/k)^{1/\gamma}} \,
\langle f_\delta, \varphi_{k,j}\rangle \varphi^{(r_1)}_k(t)\varphi^{(r_2)}_j(\tau).
\end{equation}
We note that earlier the idea of a hyperbolic cross in the problem of numerical differentiation was used in the papers
\cite{Sol_Stas_JC2020}, \cite{Sol_Stas_UMZ2022}, \cite{SSS_CMAM} (for more details about usage of hyperbolic cross in
solving the other types of ill-posed problems see \cite{PS1996}, \cite{Mileiko_Solodkii_2014},
\cite{Mileiko_Solodkii_2017_UMJ}).

The approximation properties of the method (\ref{ModVer}) will be investigated in Sections 2 and 3 while  in Section 4
it will be established that the method (\ref{ModVer}) is order-optimal in the sense of the minimal radius of the
Galerkin information.

The parameters $n$ and $\gamma$ in (\ref{ModVer}) should  be chosen depending on $\delta$, $p$, $s$, $r_1, r_2$ and $\mu$ so as to minimize  the error of the method $\mathcal{D}_{n,\gamma}^{(r_1,r_2)}$:
$$
f^{(r_1,r_2)}(t,\tau)-\mathcal{D}_{n,\gamma}^{(r_1,r_2)} f_\delta(t,\tau)
$$
\begin{equation}\label{fullError}
= \left(f^{(r_1,r_2)}(t,\tau)-\mathcal{D}_{n,\gamma}^{(r_1,r_2)}
f(t,\tau)\right)+\left(\mathcal{D}_{n,\gamma}^{(r_1,r_2)} f(t,\tau)-\mathcal{D}_{n,\gamma}^{(r_1,r_2)}
f_\delta(t,\tau)\right).
\end{equation}
For the first difference on the right-hand side of (\ref{fullError}), the representation
\begin{equation}\label{Bound_err}
 f^{(r_1,r_2)}(t,\tau)-\mathcal{D}_{n,\gamma}^{(r_1,r_2)} f(t,\tau)= \triangle_{1}(t,\tau)+\triangle_{2}(t,\tau)
\end{equation}
 holds true, where
\begin{equation}\label{Triangle_1HC}
\triangle_{1}(t,\tau)= \sum_{k=n/r_2^\gamma+1}^{\infty}\quad \sum_{j=r_2}^{\infty} \langle f,
\varphi_{kj}\rangle\varphi^{(r_1)}_k(t)\, \varphi^{(r_2)}_j(\tau),
\end{equation}
\begin{equation}\label{Triangle_2HC}
 \triangle_{2}(t,\tau)= \sum_{k=r_1}^{n/r_2^\gamma}\quad \sum_{j=(n/k)^{1/\gamma}+1}^{\infty} \langle f,
 \varphi_{kj}\rangle\varphi^{(r_1)}_k(t)\, \varphi^{(r_2)}_j(\tau) .
\end{equation}

For our calculations, we need the following formula (see Lemma 18 \cite{Mul69})
\begin{equation}\label{Muller}
 \varphi_{k}'(t) \, = 2 \, \sqrt{k+1/2}
  \mathop{{\sum}^*}\limits_{l=0}^{k-1} \sqrt{l+1/2} \, \varphi_{l}(t) ,
  \ \ \ k\in\mathbb{N} ,
 \end{equation}
where in aggregate \quad $ \mathop{{\sum}^*}\limits_{l=0}^{k-1} \sqrt{l+1/2} \, \varphi_{l}(t) \, $ \quad the summation
is extended  over only  those terms for which $k+l$ is odd.

Let us estimate the error of the method (\ref{ModVer}) in the metric of $L_2$. Upper bounds for difference
(\ref{Bound_err}) are contained in the following statement.

\begin{lemma}\label{lemma_BoundErrHC}
Let $f\in L^\mu_{s,2}$, $1\leq s< \infty$, $r_1 \geq r_2$, $\mu>2r_1-1/s+1/2$.

I)  If $r_1 = r_2$, then for $\gamma=1$ it holds true
$$
  \|f^{(r_1,r_1)}-\mathcal{D}_{n,1}^{(r_1,r_1)} f\|_{L_2}\leq c\|f\|_{s,\mu} n^{-\mu+2r_1-1/s+1/2} \ln^{3/2-1/s} n.
$$

II)  If $r_1 > r_2$, then\\
a) for any \, $\gamma \in \left[1, \frac{\mu-2r_2+1/s-1/2}{\mu-2r_1+1/s+1/2}\right) \bigcup \left(
\frac{\mu-2r_2+1/s-1/2}{\mu-2r_1+1/s+1/2}, \frac{\mu-2r_2+1/s+1/2}{\mu-2r_1+1/s+1/2}\right) \bigcup \left(
\frac{\mu-2r_2+1/s+1/2}{\mu-2r_1+1/s+1/2}, \frac{\mu-2r_2+1/s-1/2}{\mu-2r_1+1/s-1/2}\right) $ \, it holds true
$$
  \|f^{(r_1,r_2)}-\mathcal{D}_{n,\gamma}^{(r_1,r_2)} f\|_{L_2}\leq c\|f\|_{s,\mu} n^{-\mu+2r_1-1/s+1/2} ,
$$
b) for \, $\gamma = \frac{\mu-2r_2+1/s-1/2}{\mu-2r_1+1/s+1/2}$\, and \, $\gamma =
\frac{\mu-2r_2+1/s-1/2}{\mu-2r_1+1/s-1/2}$ it holds true
$$
  \|f^{(r_1,r_2)}-\mathcal{D}_{n,\gamma}^{(r_1,r_2)} f\|_{L_2}\leq c\|f\|_{s,\mu} n^{-\mu+2r_1-1/s+1/2} \ln^{1/2} n ,
$$
c) for \, $\gamma = \frac{\mu-2r_2+1/s+1/2}{\mu-2r_1+1/s+1/2}$ \, it holds true
$$
  \|f^{(r_1,r_2)}-\mathcal{D}_{n,\gamma}^{(r_1,r_2)} f\|_{L_2}\leq c\|f\|_{s,\mu} n^{-\mu+2r_1-1/s+1/2} \ln^{1-1/s} n .
$$
\end{lemma}

\textit{Proof.} Using the formula (\ref{Muller}), from (\ref{Triangle_1HC}) we have
$$
\triangle_{1}(t,\tau) = 2^{r_1+r_2} \sum_{k=n/r_2^\gamma+1}^{\infty}\quad \sum_{j=r_2}^{\infty}
\sqrt{k+1/2}\sqrt{j+1/2}\:\langle f, \varphi_{kj}\:\rangle
$$
$$
\times \mathop{{\sum}^*}\limits_{l_1=r_1-1}^{k-1}  (l_1+1/2)\: \mathop{{\sum}^*}\limits_{l_2=r_1-2}^{l_1-1}
(l_2+1/2)\ldots \mathop{{\sum}^*}\limits_{l_{r_1-1}=1}^{l_{r_1-2}-1}
(l_{r_1-1}+1/2)\mathop{{\sum}^*}\limits_{l_{r_1}=0}^{l_{r_1-1}-1} \sqrt{l_{r_1}+1/2} \,\varphi_{l_{r_1}}(t) $$
$$
 \times \mathop{{\sum}^*}\limits_{m_1=r_2-1}^{j-1} (m_1+1/2) \mathop{{\sum}^*}\limits_{m_2=r_2-2}^{m_1-1} (m_2+1/2)\ldots
 \mathop{{\sum}^*}\limits_{m_{r_2-1}=1}^{m_{r_2-2}-1}  (m_{r_2-1}+1/2) \mathop{{\sum}^*}\limits_{m_{r_2}=0}^{m_{r_2-1}-1}  \sqrt{m_{r_2}+1/2}\:
\varphi_{m_{r_2}}(\tau) .
$$
Further, we change the order of summation and get
$$
\triangle_{1}(t,\tau)= \triangle_{11}(t,\tau)+\triangle_{12}(t,\tau),
$$
where
$$
\triangle_{11}(t,\tau)= 2^{r_1+r_2} \mathop{{\sum}^*}\limits_{l_{r_1}=0}^{n/r_2^\gamma-r_1+1}  \sqrt{l_{r_1}+1/2} \:
\varphi_{l_{r_1}}(t) \mathop{{\sum}^*}\limits_{m_{r_2}=0}^{\infty}  \sqrt{m_{r_2}+1/2} \: \varphi_{m_{r_2}}(\tau)
$$
\begin{equation}\label{triangle_{11}}
\times \sum_{k=n/r_2^\gamma+1}^{\infty}\quad \sum_{j=m_{r_2}+r_2}^{\infty} \sqrt{k+1/2}\sqrt{j+1/2}\:\langle f,
\varphi_{kj}\:
 \rangle B^{r_1, r_2}_{k,j} ,
\end{equation}
$$
\triangle_{12}(t,\tau)= 2^{r_1+r_2} \mathop{{\sum}^*}\limits_{l_{r_1}=n/r_2^\gamma-r_1+2}^{\infty}  \sqrt{l_{r_1}+1/2}
\: \varphi_{l_{r_1}}(t) \mathop{{\sum}^*}\limits_{m_{r_2}=0}^{\infty}  \sqrt{m_{r_2}+1/2} \: \varphi_{m_{r_2}}(\tau)
$$
\begin{equation}\label{triangle_{12}}
 \times \sum_{k=l_{r_1}+r_1}^{\infty}\quad \sum_{j=m_{r_2}+r_2}^{\infty} \sqrt{k+1/2}\sqrt{j+1/2}\:\langle f,
\varphi_{kj}\:
 \rangle B^{r_1, r_2}_{k,j}
\end{equation}
and
$$
B^{r_1, r_2}_{k,j}:=\mathop{{\sum}^*}\limits_{l_1=l_{r_1}+r_1-1}^{k-1}
(l_1+1/2)\mathop{{\sum}^*}\limits_{l_2=l_{r_1}+r_1-2}^{l_{1}-1} (l_2+1/2) \ldots
\mathop{{\sum}^*}\limits_{l_{{r_1}-1}=l_{r_1}+1}^{l_{{r_1}-2}-1}  (l_{{r_1}-1}+1/2)
$$
$$
 \times \mathop{{\sum}^*}\limits_{m_1=m_{r_2}+r_2-1}^{j-1} (m_1+1/2) \mathop{{\sum}^*}\limits_{m_2=m_{r_2}+r_2-2}^{m_1-1} (m_2+1/2)\ldots
 \mathop{{\sum}^*}\limits_{m_{{r_2}-1}=m_{r_2}+1}^{m_{{r_2}-2}-1} (m_{{r_2}-1}+1/2)
$$
\begin{equation}\label{axular1}
 \leq c k^{2(r_1-1)} j^{2(r_2-1)} .
\end{equation}

Let's $1<s<\infty.$ At first we bound $\triangle_{1}(t,\tau)$. Using H\"older inequality and (\ref{axular1}), for
$\mu>2r_1+\frac{s-1}{s}-1/2$ we have
$$
\|\triangle_{11}\|_{L_2}^2\leq 4^{r_1+r_2} \mathop{{\sum}}\limits_{l_{r_1}=0}^{n/r_2^\gamma-r_1+1}  (l_{r_1}+1/2)
\mathop{{\sum}}\limits_{m_{r_2}=0}^{\infty}  (m_{r_2}+1/2)
$$
$$
\times \left(\sum_{k=n/r_2^\gamma+1}^{\infty}\quad \sum_{j=m_{r_2}+r_2}^{\infty} (kj)^{\mu} \:|\langle f,
\varphi_{kj}\: \rangle| \frac{B_{k,j}^{r_1, r_2}} {(kj)^{\mu-1/2}}\right)^2 \leq c \|f\|_{s,\mu}^2
\mathop{{\sum}}\limits_{l_{r_1}=0}^{n/r_2^\gamma-r_1+1}  (l_{r_1}+1/2) \mathop{{\sum}}\limits_{m_{r_2}=0}^{\infty}
(m_{r_2}+1/2)
$$
$$
 \times \left(\sum_{k=n/r_2^\gamma+1}^{\infty} k^{-(\mu-2r_1+3/2)s/(s-1)} \sum_{j=m_{r_2}+r_2}^{\infty}
 j^{-(\mu-2r_2+3/2)s/(s-1)}\right)^{2(s-1)/s}
$$
$$
\leq c  \|f\|_{s,\mu}^2 n^{-2(\mu-2r_1+3/2)+\frac{2(s-1)}{s}}  \mathop{{\sum}}\limits_{l_{r_1}=0}^{n/r_2^\gamma-r_1+1}
(l_{r_1}+1/2) \mathop{{\sum}}\limits_{m_{r_2}=0}^{\infty}  (m_{r_2}+1/2)^{-2(\mu-2r_2+3/2)+\frac{2(s-1)}{s}+1}
$$
$$
\leq c \|f\|_{s,\mu}^2 n^{-2(\mu-2r_1+1/s-1/2)} ,
$$
$$
\|\triangle_{12}\|_{L_2}^2\leq c \|f\|_{s,\mu}^2 \mathop{{\sum}}\limits_{l_{r_1}=n/r_2^\gamma-r_1+2}^{\infty}
(l_{r_1}+1/2) \mathop{{\sum}}\limits_{m_{r_2}=0}^{\infty}  (m_{r_2}+1/2)
$$
$$
\times \left(\sum_{k=l_{r_1}+r_1}^{\infty} k^{-(\mu-2r_1+3/2)s/(s-1)} \sum_{j=m_{r_2}+r_2}^{\infty}
 j^{-(\mu-2r_2+3/2)s/(s-1)}\right)^{2(s-1)/s}
$$
$$
\leq c \|f\|_{s,\mu}^2\mathop{{\sum}}\limits_{l_{r_1}=n/r_2^\gamma-r_1+2}^{\infty}
(l_{r_1}+1/2)^{-2(\mu-2r_1+3/2)+\frac{2(s-1)}{s}+1}
 \mathop{{\sum}}\limits_{m_{r_2}=0}^{\infty} (m_{r_2}+1/2)^{-2(\mu-2r_2+3/2)+\frac{2(s-1)}{s}+1}
 $$
$$
\leq c \|f\|_{s,\mu}^2n^{-2(\mu-2r_1+1/s-1/2)} .
$$
Combining the estimates for $\triangle_{11}(t,\tau)$ and $\triangle_{12}(t,\tau)$ we obtain
$$
\|\triangle_{1}\|_{L_2}\leq c \|f\|_{s,\mu} n^{-\mu+2r_1-1/s+1/2} .
$$

Next, using the formula (\ref{Muller}), from (\ref{Triangle_2HC}) we have
$$
\triangle_{2}(t,\tau)= 2^{r_1+r_2} \sum_{k=r_1}^{n/r_2^\gamma}\quad \sum_{j=(n/k)^{1/\gamma}+1}^{\infty}
\sqrt{k+1/2}\sqrt{j+1/2}\:\langle f, \varphi_{kj} \rangle
$$
$$
 \times \mathop{{\sum}^*}\limits_{l_1=r_1-1}^{k-1}  (l_1+1/2)\mathop{{\sum}^*}\limits_{l_2=r_1-2}^{l_1-1}  (l_2+1/2) \ldots
\mathop{{\sum}^*}\limits_{l_{{r_1}-1}=1}^{l_{{r_1}-2}-1}  (l_{{r_1}-1}+1/2)
\mathop{{\sum}^*}\limits_{l_{r_1}=0}^{l_{{r_1}-1}-1} \sqrt{l_{r_1}+1/2}\: \varphi_{l_{r_1}}(t)
$$
$$
 \times \mathop{{\sum}^*}\limits_{m_1={r_2}-1}^{j-1}  (m_1+1/2)\mathop{{\sum}^*}\limits_{m_2={r_2}-2}^{m_1-1}  (m_2+1/2) \ldots
 \mathop{{\sum}^*}\limits_{m_{{r_2}-1}=1}^{m_{{r_2}-2}-1}  (m_{{r_2}-1}+1/2)\mathop{{\sum}^*}\limits_{m_{r_2}=0}^{m_{{r_2}-1}-1}  \sqrt{m_{r_2}+1/2}\:
 \varphi_{m_{r_2}}(\tau) .
$$

Further, we change the order of summation and get
$$
\triangle_{2}(t,\tau)= \triangle_{21}(t,\tau)+\triangle_{22}(t,\tau) +\triangle_{23}(t,\tau) ,
$$
where
$$
\triangle_{21}(t,\tau)= 2^{r_1+r_2} \mathop{{\sum}^*}\limits_{l_{r_1}=0}^{n/r_2^\gamma-r_1}  \sqrt{l_{r_1}+1/2} \,
\varphi_{l_{r_1}}(t) \mathop{{\sum}^*}\limits_{m_{r_2}=0}^{(\frac{n}{l_{r_1}+r_1})^{1/\gamma}-r_2+1} \sqrt{m_{r_2}+1/2}
\: \varphi_{m_{r_2}}(\tau)
$$
$$
 \times \sum_{k=l_{r_1}+r_1}^{\frac{n}{(\underline{m}_{r_2}+r_2-1)^\gamma}}\quad
 \sum_{j=(n/k)^{1/\gamma}+1}^{\infty} \sqrt{k+1/2}\sqrt{j+1/2}\:\langle f, \varphi_{kj}\rangle\: B^{r_1, r_2}_{k,j} ,
$$
$$
\triangle_{22}(t,\tau)= 2^{r_1+r_2} \mathop{{\sum}^*}\limits_{l_{r_1}=0}^{n/r_2^\gamma-r_1}  \sqrt{l_{r_1}+1/2} \,
\varphi_{l_{r_1}}(t) \mathop{{\sum}^*}\limits_{m_{r_2}=2}^{(\frac{n}{l_{r_1}+r_1})^{1/\gamma}-r_2+2} \sqrt{m_{r_2}+1/2}
\: \varphi_{m_{r_2}}(\tau)
$$
$$
 \times \sum_{k=\frac{n}{(m_{r_2}+r_2-2)^\gamma}}^{n/r_2^\gamma}\quad \sum_{j=m_{r_2}+r_2}^{\infty} \sqrt{k+1/2}\sqrt{j+1/2}\:\langle f, \varphi_{kj}\:
 \rangle B^{r_1, r_2}_{k,j} ,
$$
$$
\triangle_{23}(t,\tau)= 2^{r_1+r_2} \mathop{{\sum}^*}\limits_{l_{r_1}=0}^{n/r_2^\gamma-r_1}  \sqrt{l_{r_1}+1/2} \,
\varphi_{l_{r_1}}(t) \mathop{{\sum}^*}\limits_{m_{r_2}=(\frac{n}{l_{r_1}+r_1})^{1/\gamma}-r_2+3}^{\infty}
\sqrt{m_{r_2}+1/2} \: \varphi_{m_{r_2}}(\tau)
$$
$$
 \times \sum_{k=l_{r_1}+r_1}^{n/r_2^\gamma}\quad \sum_{j=m_{r_2}+r_2}^{\infty} \sqrt{k+1/2}\sqrt{j+1/2}\:\langle f, \varphi_{kj}\:
 \rangle B^{r_1, r_2}_{k,j} .
$$
Let as before $1<s<\infty$. 
Now we estimate the norm of $\triangle_{2}$ in the case of  $r_1> r_2$ and $\gamma=1$:
$$
\|\triangle_{21}\|_{L_2}^2\leq c \|f\|_{s,\mu}^2n^{-2(\mu-2r_1+3/2)+\frac{2(s-1)}{s}}
\mathop{{\sum}}\limits_{l_{r_1}=0}^{n/r_2-r_1}  (l_{r_1}+1/2)
\mathop{{\sum}}\limits_{m_{r_2}=0}^{\frac{n}{l_{r_1}+r_1}-r_2+1} (m_{r_2}+1/2)
 $$
 $$
 \times \left(\sum_{k=l_{r_1}+r_1}^{\frac{n}{\underline{m}_{r_2}+r_2-1}} k^{2(r_1-r_2)s/(s-1)-1}\right)^{2(s-1)/s}
 \leq c\|f\|_{s,\mu}^2n^{-2(\mu-2r_1+3/2)+\frac{2(s-1)}{s}}
$$
$$
 \times \mathop{{\sum}}\limits_{l_{r_1}=0}^{n/r_2-r_1}  (l_{r_1}+1/2)
 \mathop{{\sum}}\limits_{m_{r_2}=0}^{\frac{n}{l_{r_1}+r_1}-r_2+1} (m_{r_2}+1/2)^{1-4(r_1-r_2)}
$$
$$
 \leq c\|f\|_{s,\mu}^2 n^{-2(\mu-2r_1+3/2)+\frac{2(s-1)}{s}}
\mathop{{\sum}}\limits_{l_{r_1}=0}^{n/r_2-r_1}  (l_{r_1}+1/2)
 = c  \|f\|_{s,\mu}^2n^{-2(\mu-2r_1+1/s-1/2)} ,
$$
$$
\|\triangle_{22}\|_{L_2}^2\leq c \|f\|_{s,\mu}^2 \mathop{{\sum}}\limits_{l_{r_1}=0}^{n/r_2-r_1}  (l_{r_1}+1/2)
\mathop{{\sum}}\limits_{m_{r_2}=2}^{\frac{n}{l_{r_1}+r_1}-r_2+2}  (m_{r_2}+1/2)
$$
$$
\times \left(\sum_{k=\frac{n}{m_{r_2}+r_2-2}}^{n/r_2} k^{-(\mu-2r_1+3/2)s/(s-1)} \sum_{j=m_{r_2}+r_2}^{\infty}
j^{-(\mu-2r_2+3/2)s/(s-1)}\right)^{2(s-1)/s}
$$
$$
 \leq c \|f\|_{s,\mu}^2 n^{-2(\mu-2r_1+3/2)+\frac{2(s-1)}{s}} \mathop{{\sum}}\limits_{l_{r_1}=0}^{n/r_2-r_1}  (l_{r_1}+1/2)
\mathop{{\sum}}\limits_{m_{r_2}=2}^{\frac{n}{l_{r_1}+r_1}-r_2+2}  (m_{r_2}+1/2)^{1-4(r_1-r_2)}
 $$
$$
 \leq c \|f\|_{s,\mu}^2 n^{-2(\mu-2r_1+1/s-1/2)} ,
$$
$$
\|\triangle_{23}\|_{L_2}^2\leq c \|f\|_{s,\mu}^2 \mathop{{\sum}}\limits_{l_{r_1}=0}^{n/r_2-r_1}  (l_{r_1}+1/2)
\mathop{{\sum}}\limits_{m_{r_2}=\frac{n}{l_{r_1}+r_1}-r_2+3}^{\infty}  (m_{r_2}+1/2)
$$
$$
\times \left(\sum_{k=l_{r_1}+r_1}^{n/r_2} k^{-(\mu-2r_1+3/2)s/(s-1)} \sum_{j=m_{r_2}+r_2}^{\infty}
j^{-(\mu-2r_2+3/2)s/(s-1)}\right)^{2(s-1)/s}
$$
$$
 \leq c \|f\|_{s,\mu}^2 \mathop{{\sum}}\limits_{l_{r_1}=0}^{n/r_2-r_1}  (l_{r_1}+1/2)^{1-2(\mu-2r_1+3/2)+\frac{2(s-1)}{s}}
\mathop{{\sum}}\limits_{m_{r_2}=\frac{n}{l_{r_1}+r_1}-r_2+3}^{\infty}
(m_{r_2}+1/2)^{1-2(\mu-2r_2+3/2)+\frac{2(s-1)}{s}}
$$
$$
 \leq c \|f\|_{s,\mu}^2 n^{-2(\mu-2r_2+3/2)+\frac{2(s-1)}{s}+2} \mathop{{\sum}}\limits_{l_{r_1}=0}^{n/r_2-r_1}
 (l_{r_1}+1/2)^{4(r_1-r_2)-1}
$$
$$
 \leq c \|f\|_{s,\mu}^2 n^{-2(\mu-2r_1+1/s-1/2)} .
$$
Combining the estimates for $\triangle_{21}(t,\tau)$, $\triangle_{22}(t,\tau)$ and $\triangle_{23}(t,\tau)$ we obtain
$$
\|\triangle_{2}\|_{L_2}\leq c\, \|f\|_{s,\mu}n^{-\mu+2r_1-1/s+1/2} .
$$
To estimate $\|\triangle_{2}\|_{L_2}$ for $r_1=r_2$ or $\gamma>1$, it is enough to repeat an analogous line of
reasoning. The combination of (\ref{Bound_err}) and bounds for the norms of $\triangle_{1}$, $\triangle_{2}$ makes it
possible to establish the desired inequality in the case $1<s<\infty$.

For $s=1$ the assertion of Lemma is proved similarly. \vspace{0.1in}

\vskip -4mm

  ${}$ \ \ \ \ \ \ \ \ \ \ \ \ \ \ \ \ \ \ \ \ \ \ \ \ \ \ \ \ \ \ \ \ \ \ \ \ \ \
  \ \ \ \ \ \ \ \ \ \ \ \ \ \ \ \ \ \ \ \ \ \
  \ \ \ \ \ \ \ \ \ \ \ \ \ \ \ \ \ \ \ \ \ \
  \ \ \ \ \ \ \ \ \ \ \ \ \ \ \ \ \ \ \ \ \ \
  \ \ \ \ \ \ \ \ \ \ \ \ \ \ \ \ \ \ \ \ \ \
  \ \ \ \ \ \ \ \ \ \ \ \ \ \ \ \ \ \ \ \ \ \
  \ \ $\Box$

The following statement contains estimates for the second difference from the right-hand side of (\ref{fullError}) in
the metric of $L_2$.

\begin{lemma}\label{lemma_BoundPertHC}
Let $f\in L_2(Q)$ and let the condition (\ref{perturbation}) be satisfied.

I)  If $r_1 = r_2$, then\\
a) for $\gamma=1$ it holds true
 $$
\|\mathcal{D}^{(r_1,r_1)}_{n,1} f - \mathcal{D}^{(r_1,r_1)}_{n,1} f_\delta\|_{L_2} \leq c\delta n^{2r_1-1/p+1/2}
\ln^{3/2-1/p} n ,
 $$
b) for any \, $1 < \gamma < \frac{2r_1-1/2+(p-1)/p}{2r_1-3/2+(p-1)/p}$ \, it holds true
 $$
\|\mathcal{D}^{(r_1,r_1)}_{n,\gamma} f - \mathcal{D}^{(r_1,r_1)}_{n,\gamma} f_\delta\|_{L_2} \leq c\delta
n^{2r_1-1/p+1/2} ,
 $$
c) for \, $\gamma = \frac{2r_1-1/2+(p-1)/p}{2r_1-3/2+(p-1)/p}$ \, it holds true
 $$
\|\mathcal{D}^{(r_1,r_1)}_{n,\gamma} f - \mathcal{D}^{(r_1,r_1)}_{n,\gamma} f_\delta\|_{L_2} \leq c\delta
n^{2r_1-1/p+1/2} \ln^{1/2} n .
 $$

II)  If $r_1 > r_2$, then for any $\gamma\geq 1$ it holds true
 $$
\|\mathcal{D}^{(r_1,r_2)}_{n,\gamma} f - \mathcal{D}^{(r_1,r_2)}_{n,\gamma} f_\delta\|_{L_2} \leq c\delta
n^{2r_1-1/p+1/2} .
 $$
\end{lemma}

\textit{Proof.} Let's write down the representation
$$
\mathcal{D}_{n,\gamma}^{(r_1,r_2)} f(t,\tau) - \mathcal{D}_{n,\gamma}^{(r_1,r_2)} f_\delta(t,\tau) =
\sum_{k=r_1}^{n/r_2^\gamma}\, \sum_{j=r_2}^{(n/k)^{1/\gamma}} \, \langle f - f_\delta, \varphi_{k,j}\rangle
\varphi^{(r_1)}_k(t)\varphi^{(r_2)}_j(\tau).
$$
Using the formula (\ref{Muller}), we get
$$\mathcal{D}_{n,\gamma}^{(r_1,r_2)} f(t,\tau) - \mathcal{D}_{n,\gamma}^{(r_1,r_2)} f_\delta(t,\tau)= 2^{r_1+r_2} \sum_{k=r_1}^{n/r_2^\gamma}\,
\sum_{j=r_2}^{(n/k)^{1/\gamma}} \sqrt{k+1/2}\sqrt{j+1/2}\:\langle f-f_\delta, \varphi_{kj}\: \rangle
$$
$$
\times \mathop{{\sum}^*}\limits_{l_1=r_1-1}^{k-1}  (l_1+1/2)\: \mathop{{\sum}^*}\limits_{l_2=r_1-2}^{l_1-1} (l_2+1/2)\:
\ldots \mathop{{\sum}^*}\limits_{l_{{r_1}-1}=1}^{l_{{r_1}-2}-1} (l_{{r_1}-1}+1/2)\:
 \mathop{{\sum}^*}\limits_{l_{r_1}=0}^{l_{{r_1}-1}-1} \sqrt{l_{r_1}+1/2}\: \varphi_{l_{r_1}}(t)
 $$
 $$
\times \mathop{{\sum}^*}\limits_{m_1={r_2}-1}^{j-1}  (m_1+1/2)\: \mathop{{\sum}^*}\limits_{m_2={r_2}-2}^{m_1-1}
(m_2+1/2)\: \ldots \mathop{{\sum}^*}\limits_{m_{{r_2}-1}=1}^{m_{{r_2}-2}-1} (m_{{r_2}-1}+1/2)\:
 \mathop{{\sum}^*}\limits_{m_{r_2}=0}^{m_{{r_2}-1}-1} \sqrt{m_{r_2}+1/2}\: \varphi_{m_{r_2}}(\tau) .
$$
Further, we change the order of summation and get
$$
\mathcal{D}_{n,\gamma}^{(r_1,r_2)} f(t,\tau) - \mathcal{D}_{n,\gamma}^{(r_1,r_2)} f_\delta(t,\tau)= 2^{r_1+r_2} \,
\mathop{{\sum}^*}\limits_{l_{r_1}=0}^{n/r_2^\gamma-r_1} \sqrt{l_{r_1}+1/2}\, \varphi_{l_{r_1}}(t)
\mathop{{\sum}^*}\limits_{m_{r_2}=0}^{(\frac{n}{l_{r_1}+r_1})^{1/\gamma}-r_2}  \sqrt{m_{r_2}+1/2} \:
\varphi_{m_{r_2}}(\tau)
$$
$$
\times \sum_{k=l_{r_1}+r_1}^{\frac{n}{(m_{r_2}+r_2)^{\gamma}}}\quad \sum_{j=m_{r_2}+r_2}^{(\frac{n}{k})^{1/\gamma}}
\sqrt{k+1/2}\sqrt{j+1/2}\, \langle f-f_\delta, \varphi_{kj}\: \rangle B^{r_1, r_2}_{k,j} .
$$
Let $1<p<\infty$ first. Then, using the H\"older inequality and the estimate (\ref{axular1}), we find
$$
\|\mathcal{D}_{n,\gamma}^{(r_1,r_2)} f - \mathcal{D}_{n,\gamma}^{(r_1,r_2)} f_\delta\|_{L_2}^2 \leq 4^{r_1+r_2}\,
\mathop{{\sum}}\limits_{l_{r_1}=0}^{n/r_2^\gamma-r_1} (l_{r_1}+1/2)\,
\mathop{{\sum}}\limits_{m_{r_2}=0}^{(\frac{n}{l_{r_1}+r_1})^{1/\gamma}-r_2}  (m_{r_2}+1/2)
$$
$$
\times \left(\sum_{k=l_{r_1}+r_1}^{\frac{n}{(m_{r_2}+r_2)^{\gamma}}}\quad
\sum_{j=m_{r_2}+r_2}^{(\frac{n}{k})^{1/\gamma}} \sqrt{k+1/2}\sqrt{j+1/2}\, |\langle f-f_\delta, \varphi_{kj}\: \rangle
| B^{r_1, r_2}_{k,j}\right)^{2}
$$
$$
\leq  c \delta^2 \mathop{{\sum}}\limits_{l_{r_1}=0}^{n/r_2^\gamma-r_1} (l_{r_1}+1/2)\,
\mathop{{\sum}}\limits_{m_{r_2}=0}^{(\frac{n}{l_{r_1}+r_1})^{1/\gamma}-r_2}  (m_{r_2}+1/2)
$$
$$ \times
\left(\sum_{k=l_{r_1}+r_1}^{\frac{n}{(m_{r_2}+r_2)^{\gamma}}} k^{(2r_1-3/2)p/(p-1)}
\sum_{j=m_{r_2}+r_2}^{(\frac{n}{k})^{1/\gamma}} j^{(2r_2-3/2)p/(p-1)}\right)^{2(p-1)/p}
$$
$$\leq
 c \delta^2 \mathop{{\sum}}\limits_{l_{r_1}=0}^{n/r_2^\gamma-r_1} (l_{r_1}+1/2)\,
\mathop{{\sum}}\limits_{m_{r_2}=0}^{(\frac{n}{l_{r_1}+r_1})^{1/\gamma}-r_2}  (m_{r_2}+1/2)
$$
\begin{equation} \label{Sec_dif_L2}
\times \left(n^{\frac{(2r_2-3/2)p}{\gamma(p-1)}+1/\gamma} \sum_{k=l_{r_1}+r_1}^{\frac{n}{(m_{r_2}+r_2)^{\gamma}}}
k^{\frac{(2r_1-3/2)p}{p-1}-\frac{(2r_2-3/2)p}{\gamma(p-1)}-1/\gamma}\right)^{2(p-1)/p} .
\end{equation}
 Whence it follows for $r_1> r_2$ and any $\gamma\geq 1$ that
$$
\|\mathcal{D}_{n,\gamma}^{(r_1,r_2)} f - \mathcal{D}_{n,\gamma}^{(r_1,r_2)} f_\delta\|_{L_2}^2 \leq c \delta^2
n^{4r_1-3+\frac{2(p-1)}{p}} \, \mathop{{\sum}}\limits_{l_{r_1}=0}^{n/r_2^\gamma-r_1} (l_{r_1}+1/2) = c \delta^2
n^{4r_1-2/p+1} .
$$
which was required to prove.

Let's now $r_1 = r_2$ and any $\gamma=1$. Then from (\ref{Sec_dif_L2}) we obtain
$$
\|\mathcal{D}_{n,1}^{(r_1,r_1)} f - \mathcal{D}_{n,1}^{(r_1,r_1)} f_\delta\|_{L_2}^2 \leq c\delta^2
n^{4r_1-3+\frac{2(p-1)}{p}} \, \ln^{\frac{2(p-1)}{p}} n\, \mathop{{\sum}}\limits_{l_{r_1}=0}^{n/r_1-r_1} (l_{r_1}+1/2)
\mathop{{\sum}}\limits_{m_{r_1}=0}^{\frac{n}{l_{r_1}+r_1}-r_1}  (m_{r_1}+1/2)
$$
$$
\leq c \delta^2 n^{4r_1-1+\frac{2(p-1)}{p}} \, \ln^{\frac{2(p-1)}{p}} n\,
\mathop{{\sum}}\limits_{l_{r_1}=0}^{n/r_2-r_1} \frac{1}{l_{r_1}+1/2} = c \delta^2 n^{4r_1-2/p+1}\, \ln^{3-2/p} n .
$$

To estimate $\|\mathcal{D}_{n,\gamma}^{(r_1,r_2)} f - \mathcal{D}_{n,\gamma}^{(r_1,r_2)} f_\delta\|_{L_2}$ for $\gamma
>1$ and also for $p=1$, $p=\infty$, it is enough to repeat an analogous line of reasoning.  \vspace{0.1in}

\vskip -4mm

  ${}$ \ \ \ \ \ \ \ \ \ \ \ \ \ \ \ \ \ \ \ \ \ \ \ \ \ \ \ \ \ \ \ \ \ \ \ \ \ \
  \ \ \ \ \ \ \ \ \ \ \ \ \ \ \ \ \ \ \ \ \ \
  \ \ \ \ \ \ \ \ \ \ \ \ \ \ \ \ \ \ \ \ \ \
  \ \ \ \ \ \ \ \ \ \ \ \ \ \ \ \ \ \ \ \ \ \
  \ \ \ \ \ \ \ \ \ \ \ \ \ \ \ \ \ \ \ \ \ \
  \ \ \ \ \ \ \ \ \ \ \ \ \ \ \ \ \ \ \ \ \ \
  \ \ $\Box$

The combination of Lemmas \ref{lemma_BoundErrHC} and \ref{lemma_BoundPertHC} gives
\begin{theorem} \label{Th1}
Let $f\in L^\mu_{s,2}$, $1\leq s< \infty$, $r_1 \geq r_2$, $\mu>2r_1+1/2-1/s$, and let the condition
(\ref{perturbation}) be satisfied.

I)  If $r_1 = r_2$, then for $n\asymp \left(\delta/ \ln^{1/p-1/s} \frac{1}{\delta}\right)^{-\frac{1}{\mu-1/p+1/s}}$ and
$\gamma=1$ it holds true
  $$
   \|f^{(r_1,r_1)} - \mathcal{D}^{(r_1,r_1)}_{n,1} f_\delta\|_{L_2} \leq c \left(\delta \ln^{1/s-1/p}
\frac{1}{\delta}\right)^{\frac{\mu-2r_1+1/s-1/2}{\mu-1/p+1/s}} \ln^{3/2-1/s} \frac{1}{\delta} .
  $$

II)  If $r_1 > r_2$, then\\
a) for $n\asymp \delta^{-\frac{1}{\mu-1/p+1/s}}$ and any \, $ \gamma \in \left[1,
\frac{\mu-2r_2+1/s-1/2}{\mu-2r_1+1/s+1/2}\right) \bigcup \left( \frac{\mu-2r_2+1/s-1/2}{\mu-2r_1+1/s+1/2},
\frac{\mu-2r_2+1/s+1/2}{\mu-2r_1+1/s+1/2}\right) $\\
$
 \bigcup \left(
\frac{\mu-2r_2+1/s+1/2}{\mu-2r_1+1/s+1/2}, \frac{\mu-2r_2+1/s-1/2}{\mu-2r_1+1/s-1/2}\right) $ \, it holds true
$$
  \|f^{(r_1,r_2)}-\mathcal{D}_{n,\gamma}^{(r_1,r_2)} f_\delta\|_{L_2}\leq c \delta
^{\frac{\mu-2r_1+1/s-1/2}{\mu-1/p+1/s}} ,
$$
b) for $n\asymp \delta^{-\frac{1}{\mu-1/p+1/s}} \ln^{\frac{1}{2(\mu-1/p+1/s)}} \frac{1}{\delta}$ and $\, \gamma =
\frac{\mu-2r_2+1/s-1/2}{\mu-2r_1+1/s+1/2}, \quad \gamma = \frac{\mu-2r_2+1/s-1/2}{\mu-2r_1+1/s-1/2}\, $\\ it holds true
$$
  \|f^{(r_1,r_2)}-\mathcal{D}_{n,\gamma}^{(r_1,r_2)} f_\delta\|_{L_2}\leq c \left(\delta/ \ln^{1/2}
\frac{1}{\delta}\right)^{\frac{\mu-2r_1+1/s-1/2}{\mu-1/p+1/s}} \ln^{1/2} \frac{1}{\delta}  ,
$$
c) for $n\asymp \left(\delta/ \ln^{1-1/s} \frac{1}{\delta}\right)^{-\frac{1}{\mu-1/p+1/s}}$ and $\, \gamma =
\frac{\mu-2r_2+1/s+1/2}{\mu-2r_1+1/s+1/2}\, $ it holds true
$$
  \|f^{(r_1,r_2)}-\mathcal{D}_{n,\gamma}^{(r_1,r_2)} f_\delta\|_{L_2}\leq c \left(\delta/ \ln^{1-1/s}
\frac{1}{\delta}\right)^{\frac{\mu-2r_1+1/s-1/2}{\mu-1/p+1/s}} \ln^{1-1/s} \frac{1}{\delta}  .
$$
\end{theorem}

\vskip 2mm
\section{Truncation method. Error estimate in the metric of $C$}

Now we have to bound the error of (\ref{ModVer}) in the metric of $C$. Upper estimates for the norm of difference
(\ref{Bound_err}) are contained in the following statement.

\begin{lemma}\label{lemma_BoundErrHCC}
Let $f\in L^\mu_{s,2}$, $1\leq s< \infty$, $r_1 \geq r_2$, $\mu>2r_1+3/2-1/s$.

I)  If $r_1 = r_2$, then for $\gamma=1$ it holds true
$$
  \|f^{(r_1,r_1)}-\mathcal{D}_{n,1}^{(r_1,r_1)} f\|_{C}\leq c\|f\|_{s,\mu} n^{-\mu+2r_1-1/s+3/2} \ln^{2-1/s} n.
$$

II)  If $r_1 = r_2+1$, then\\
a) for any \, $\gamma \in \left(1, \frac{\mu-2r_1+1/s+5/2}{\mu-2r_1+1/s+1/2}\right) \bigcup \left(
\frac{\mu-2r_1+1/s+5/2}{\mu-2r_1+1/s+1/2}, \frac{\mu-2r_1+1/s+1/2}{\mu-2r_1+1/s-3/2}\right) $\\ it holds true
$$
  \|f^{(r_1,r_1-1)}-\mathcal{D}_{n,\gamma}^{(r_1,r_1-1)} f\|_{C}\leq c\|f\|_{s,\mu} n^{-\mu+2r_1-1/s+3/2} ,
$$
b) for \, $\gamma = 1$ \, and\,  $\gamma = \frac{\mu-2r_1+1/s+1/2}{\mu-2r_1+1/s-3/2} $ it holds true
$$
  \|f^{(r_1,r_1-1)}-\mathcal{D}_{n,\gamma}^{(r_1,r_1-1)} f\|_{C} \leq c\|f\|_{s,\mu} n^{-\mu+2r_1-1/s+3/2} \ln n ,
$$
c) for \, $\gamma = \frac{\mu-2r_1+1/s+5/2}{\mu-2r_1+1/s+1/2} $\, it holds true
$$
  \|f^{(r_1,r_1-1)}-\mathcal{D}_{n,\gamma}^{(r_1,r_1-1)} f\|_{C} \leq c\|f\|_{s,\mu} n^{-\mu+2r_1-1/s+3/2} \ln^{1-1/s} n .
$$

III)  If $r_1 > r_2+1$, then\\
a) for any \, $\gamma \in \left[1, \frac{\mu-2r_2+1/s-3/2}{\mu-2r_1+1/s+1/2}\right) \bigcup
\left( \frac{\mu-2r_2+1/s-3/2}{\mu-2r_1+1/s+1/2}, \frac{\mu-2r_2+1/s+1/2}{\mu-2r_1+1/s+1/2}\right) $\\
$
 \bigcup \left(
\frac{\mu-2r_2+1/s+1/2}{\mu-2r_1+1/s+1/2}, \frac{\mu-2r_2+1/s-3/2}{\mu-2r_1+1/s-3/2}\right) $\, it holds true
$$
  \|f^{(r_1,r_2)}-\mathcal{D}_{n,\gamma}^{(r_1,r_2)} f\|_{C}\leq c\|f\|_{s,\mu} n^{-\mu+2r_1-1/s+3/2} ,
$$
b) for \, $\gamma = \frac{\mu-2r_2+1/s-3/2}{\mu-2r_1+1/s+1/2}$\,  and \, $\gamma =
\frac{\mu-2r_2+1/s-3/2}{\mu-2r_1+1/s-3/2} $ \,it holds true
$$
  \|f^{(r_1,r_2)}-\mathcal{D}_{n,\gamma}^{(r_1,r_2)} f\|_{C}\leq c\|f\|_{s,\mu} n^{-\mu+2r_1-1/s+3/2} \ln n ,
$$
c) for \, $\gamma = \frac{\mu-2r_2+1/s+1/2}{\mu-2r_1+1/s+1/2} $\, it holds true
$$
  \|f^{(r_1,r_2)}-\mathcal{D}_{n,\gamma}^{(r_1,r_2)} f\|_{C}\leq c\|f\|_{s,\mu} n^{-\mu+2r_1-1/s+3/2} \ln^{1-1/s} n .
$$
\end{lemma}

\textit{Proof.} Let's start with the case $1< s< \infty$. Using (\ref{axular1}), from (\ref{triangle_{11}}) and
(\ref{triangle_{12}})  we get
$$
\|\triangle_{11}\|_{C} \leq c \mathop{{\sum}}\limits_{l_{r_1}=0}^{n/r_2^\gamma-r_1+1}  (l_{r_1}+1/2)
\mathop{{\sum}}\limits_{m_{r_2}=0}^{\infty}  (m_{r_2}+1/2)
$$
$$
\times \sum_{k=n/r_2^\gamma+1}^{\infty}\quad \sum_{j=m_{r_2}+r_2}^{\infty} (kj)^{\mu} \:|\langle f, \varphi_{kj}\:
\rangle| \frac{B_{k,j}^{r_1, r_2}} {(kj)^{\mu-1/2}} \leq c \|f\|_{s,\mu}
\mathop{{\sum}}\limits_{l_{r_1}=0}^{n/r_2^\gamma-r_1+1}  (l_{r_1}+1/2) \mathop{{\sum}}\limits_{m_{r_2}=0}^{\infty}
(m_{r_2}+1/2)
$$
$$
 \times \left(\sum_{k=n/r_2^\gamma+1}^{\infty} k^{-(\mu-2r_1+3/2)s/(s-1)} \sum_{j=m_{r_2}+r_2}^{\infty}
 j^{-(\mu-2r_2+3/2)s/(s-1)}\right)^{(s-1)/s}
$$
$$
\leq c  \|f\|_{s,\mu} n^{-(\mu-2r_1+3/2)+\frac{s-1}{s}}  \mathop{{\sum}}\limits_{l_{r_1}=0}^{n/r_2^\gamma-r_1+1}
(l_{r_1}+1/2) \mathop{{\sum}}\limits_{m_{r_2}=0}^{\infty}  (m_{r_2}+1/2)^{1-(\mu-2r_2+3/2)+\frac{s-1}{s}}
$$
$$
\leq c \|f\|_{s,\mu} n^{-\mu+2r_1-1/s+3/2} ,
$$
$$
\|\triangle_{12}\|_{C} \leq c \|f\|_{s,\mu}  \mathop{{\sum}}\limits_{l_{r_1}=n/r_2^\gamma-r_1+2}^{\infty} (l_{r_1}+1/2)
\mathop{{\sum}}\limits_{m_{r_2}=0}^{\infty}  (m_{r_2}+1/2)
$$
$$
\times \left(\sum_{k=l_{r_1}+r_1}^{\infty} k^{-(\mu-2r_1+3/2)s/(s-1)} \sum_{j=m_{r_2}+r_2}^{\infty}
 j^{-(\mu-2r_2+3/2)s/(s-1)}\right)^{(s-1)/s}
$$
$$
\leq c \|f\|_{s,\mu} \mathop{{\sum}}\limits_{l_{r_1}=n/r_2^\gamma-r_1+2}^{\infty}
(l_{r_1}+1/2)^{1-(\mu-2r_1+3/2)+\frac{s-1}{s}}
 \mathop{{\sum}}\limits_{m_{r_2}=0}^{\infty} (m_{r_2}+1/2)^{1-(\mu-2r_2+3/2)+\frac{s-1}{s}}
 $$
$$
\leq c \|f\|_{s,\mu} n^{-\mu+2r_1-1/s+3/2} .
$$
Thus, we get
$$
\|\triangle_{1}\|_{C} \leq c \|f\|_{s,\mu}  n^{-\mu+2r_1-1/s+3/2} .
$$
Further, we estimate the norm of $\triangle_{2}$ in the case of  $r_1> r_2$ and $\gamma=1$:
$$
\|\triangle_{21}\|_{C} \leq c \|f\|_{s,\mu} n^{-\mu+2r_2-3/2+\frac{s-1}{s}}
\mathop{{\sum}}\limits_{l_{r_1}=0}^{n/r_2-r_1}  (l_{r_1}+1/2)
\mathop{{\sum}}\limits_{m_{r_2}=0}^{\frac{n}{l_{r_1}+r_1}-r_2+1} (m_{r_2}+1/2)
 $$
 $$
 \times \left(\sum_{k=l_{r_1}+r_1}^{\frac{n}{\underline{m}_{r_2}+r_2-1}} k^{2(r_1-r_2)s/(s-1)-1}\right)^{(s-1)/s}
 \leq c\|f\|_{s,\mu} n^{-\mu+2r_1-3/2+\frac{s-1}{s}}
$$
$$
 \times \mathop{{\sum}}\limits_{l_{r_1}=0}^{n/r_2-r_1}  (l_{r_1}+1/2)
 \mathop{{\sum}}\limits_{m_{r_2}=0}^{\frac{n}{l_{r_1}+r_1}-r_2+1} (m_{r_2}+1/2)^{1-2(r_1-r_2)}
$$
$$
 \leq c\|f\|_{s,\mu} \left\{
\begin{array}{cl}
n^{-\mu+2r_1-1/s+3/2}\, \ln n , & r_1=r_2+1, \\
n^{-\mu+2r_1-1/s+3/2} , & r_1>r_2+1,
\end{array}
\right.
$$
$$
\|\triangle_{22}\|_{C} \leq c \|f\|_{s,\mu} \mathop{{\sum}}\limits_{l_{r_1}=0}^{n/r_2-r_1}  (l_{r_1}+1/2)
\mathop{{\sum}}\limits_{m_{r_2}=2}^{\frac{n}{l_{r_1}+r_1}-r_2+2}  (m_{r_2}+1/2)
$$
$$
\times \left(\sum_{k=\frac{n}{m_{r_2}+r_2-2}}^{n/r_2} k^{-(\mu-2r_1+3/2)s/(s-1)} \sum_{j=m_{r_2}+r_2}^{\infty}
j^{-(\mu-2r_2+3/2)s/(s-1)}\right)^{(s-1)/s}
$$
$$
 \leq c \|f\|_{s,\mu} n^{-\mu+2r_1-3/2+\frac{s-1}{s}} \mathop{{\sum}}\limits_{l_{r_1}=0}^{n/r_2-r_1}  (l_{r_1}+1/2)
\mathop{{\sum}}\limits_{m_{r_2}=2}^{\frac{n}{l_{r_1}+r_1}-r_2+2}  (m_{r_2}+1/2)^{1-2(r_1-r_2)}
$$
$$
 \leq c\|f\|_{s,\mu} \left\{
\begin{array}{cl}
n^{-\mu+2r_1-1/s+3/2}\, \ln n , & r_1=r_2+1, \\
n^{-\mu+2r_1-1/s+3/2} , & r_1>r_2+1,
\end{array}
\right.
$$
$$
\|\triangle_{23}\|_{C} \leq c \|f\|_{s,\mu} \mathop{{\sum}}\limits_{l_{r_1}=0}^{n/r_2-r_1}  (l_{r_1}+1/2)
\mathop{{\sum}}\limits_{m_{r_2}=\frac{n}{l_{r_1}+r_1}-r_2+3}^{\infty}  (m_{r_2}+1/2)
$$
$$
\times \left(\sum_{k=l_{r_1}+r_1}^{n/r_2} k^{-(\mu-2r_1+3/2)s/(s-1)} \sum_{j=m_{r_2}+r_2}^{\infty}
j^{-(\mu-2r_2+3/2)s/(s-1)}\right)^{(s-1)/s}
$$
$$
 \leq c \|f\|_{s,\mu} \mathop{{\sum}}\limits_{l_{r_1}=0}^{n/r_2-r_1}  (l_{r_1}+1/2)^{1-(\mu-2r_1+3/2)+\frac{s-1}{s}}
\mathop{{\sum}}\limits_{m_{r_2}=\frac{n}{l_{r_1}+r_1}-r_2+3}^{\infty} (m_{r_2}+1/2)^{1-(\mu-2r_2+3/2)+\frac{s-1}{s}}
$$
$$
 \leq c \|f\|_{s,\mu} n^{-(\mu-2r_2+3/2)+\frac{s-1}{s}+2} \mathop{{\sum}}\limits_{l_{r_1}=0}^{n/r_2-r_1}
 (l_{r_1}+1/2)^{2(r_1-r_2)-1}
$$
$$
 \leq c \|f\|_{s,\mu} n^{-\mu+2r_1-1/s+3/2} .
$$
Combining the estimates for $\triangle_{21}(t,\tau)$, $\triangle_{22}(t,\tau)$ and $\triangle_{23}(t,\tau)$ we obtain
$$
\|\triangle_{2}\|_{C}  \leq c\|f\|_{s,\mu} \left\{
\begin{array}{cl}
n^{-\mu+2r_1-1/s+3/2}\, \ln n , & r_1=r_2+1, \\
n^{-\mu+2r_1-1/s+3/2} , & r_1>r_2+1,
\end{array}
\right. .
$$
To estimate $\|\triangle_{2}\|_{C}$ for $r_1=r_2$ or $\gamma>1$, it is enough to repeat an analogous line of reasoning.
The combination of (\ref{Bound_err}) and bounds for the norms of $\triangle_{1}$, $\triangle_{2}$ makes it possible to
establish the desired inequality in the case of $1<s<\infty$.

For $s=1$ the assertion of Lemma is proved similarly. \vspace{0.1in}

\vskip -4mm

  ${}$ \ \ \ \ \ \ \ \ \ \ \ \ \ \ \ \ \ \ \ \ \ \ \ \ \ \ \ \ \ \ \ \ \ \ \ \ \ \
  \ \ \ \ \ \ \ \ \ \ \ \ \ \ \ \ \ \ \ \ \ \
  \ \ \ \ \ \ \ \ \ \ \ \ \ \ \ \ \ \ \ \ \ \
  \ \ \ \ \ \ \ \ \ \ \ \ \ \ \ \ \ \ \ \ \ \
  \ \ \ \ \ \ \ \ \ \ \ \ \ \ \ \ \ \ \ \ \ \
  \ \ \ \ \ \ \ \ \ \ \ \ \ \ \ \ \ \ \ \ \ \
  \ \ $\Box$

The following statement contains upper estimates for the second difference from the right-hand side of
(\ref{fullError}) in the metric of $C$.
\begin{lemma}\label{lemma_BoundPertHCC}
Let $f\in L_2(Q)$ and let the condition (\ref{perturbation}) be satisfied.

I)  If $r_1 = r_2$, then\\
a) for $\gamma=1$ it holds true
 $$
\|\mathcal{D}^{(r_1,r_1)}_{n,1} f - \mathcal{D}^{(r_1,r_1)}_{n,1} f_\delta\|_{C} \leq c\delta n^{2r_1-1/p+3/2}
\ln^{2-1/p} n ,
 $$
b) for any \, $\quad 1 < \gamma < \frac{2r_1+1/2+(p-1)/p}{2r_1-3/2+(p-1)/p}$ \, it holds true
 $$
\|\mathcal{D}^{(r_1,r_1)}_{n,\gamma} f - \mathcal{D}^{(r_1,r_1)}_{n,\gamma} f_\delta\|_{C} \leq c\delta
n^{2r_1-1/p+3/2} ,
 $$
c) for \, $\gamma = \frac{2r_1+1/2+(p-1)/p}{2r_1-3/2+(p-1)/p}$ \, it holds true
 $$
\|\mathcal{D}^{(r_1,r_1)}_{n,\gamma} f - \mathcal{D}^{(r_1,r_1)}_{n,\gamma} f_\delta\|_{C} \leq c\delta
n^{2r_1-1/p+3/2} \ln n .
 $$

II)  If $r_1 = r_2+1$, then\\
a) for $\gamma= 1$ it holds true
 $$
\|\mathcal{D}^{(r_1,r_1-1)}_{n,\gamma} f - \mathcal{D}^{(r_1,r_1-1)}_{n,\gamma} f_\delta\|_{C} \leq c\delta
n^{2r_1-1/p+3/2} \ln n ,
 $$
b) for any \, $\gamma >1$ \, it holds true
 $$
\|\mathcal{D}^{(r_1,r_1-1)}_{n,\gamma} f - \mathcal{D}^{(r_1,r_1-1)}_{n,\gamma} f_\delta\|_{C} \leq c\delta
n^{2r_1-1/p+3/2} .
 $$

III)  If $r_1 > r_2+1$, then for any $\gamma\geq 1$ it holds true
 $$
\|\mathcal{D}^{(r_1,r_2)}_{n,\gamma} f - \mathcal{D}^{(r_1,r_2)}_{n,\gamma} f_\delta\|_{C} \leq c\delta
n^{2r_1-1/p+3/2} .
 $$
\end{lemma}
\textit{Proof.} Let $1< p<\infty$ first. Then, using the H\"older inequality and the estimate (\ref{axular1}), we find
$$
\|\mathcal{D}_{n,\gamma}^{(r_1,r_2)} f - \mathcal{D}_{n,\gamma}^{(r_1,r_2)} f_\delta\|_{C} \leq 2^{r_1+r_2}\,
\mathop{{\sum}}\limits_{l_{r_1}=0}^{n/r_2^\gamma-r_1} (l_{r_1}+1/2)\,
\mathop{{\sum}}\limits_{m_{r_2}=0}^{(\frac{n}{l_{r_1}+r_1})^{1/\gamma}-r_2}  (m_{r_2}+1/2)
$$
$$
\times \sum_{k=l_{r_1}+r_1}^{\frac{n}{(m_{r_2}+r_2)^{\gamma}}}\quad \sum_{j=m_{r_2}+r_2}^{(\frac{n}{k})^{1/\gamma}}
\sqrt{k+1/2}\sqrt{j+1/2}\, |\langle f-f_\delta, \varphi_{kj}\: \rangle | B^{r_1, r_2}_{k,j}
$$
$$
\leq  c \delta \mathop{{\sum}}\limits_{l_{r_1}=0}^{n/r_2^\gamma-r_1} (l_{r_1}+1/2)\,
\mathop{{\sum}}\limits_{m_{r_2}=0}^{(\frac{n}{l_{r_1}+r_1})^{1/\gamma}-r_2}  (m_{r_2}+1/2)
$$
$$ \times
\left(\sum_{k=l_{r_1}+r_1}^{\frac{n}{(m_{r_2}+r_2)^{\gamma}}} k^{(2r_1-3/2)p/(p-1)}
\sum_{j=m_{r_2}+r_2}^{(\frac{n}{k})^{1/\gamma}} j^{(2r_2-3/2)p/(p-1)}\right)^{(p-1)/p}
$$
$$\leq
 c \delta n^{(2r_2-3/2)/\gamma+\frac{p-1}{\gamma p}} \, \mathop{{\sum}}\limits_{l_{r_1}=0}^{n/r_2^\gamma-r_1} (l_{r_1}+1/2)\,
\mathop{{\sum}}\limits_{m_{r_2}=0}^{(\frac{n}{l_{r_1}+r_1})^{1/\gamma}-r_2}  (m_{r_2}+1/2)
$$
$$
\times \left(\sum_{k=l_{r_1}+r_1}^{\frac{n}{(m_{r_2}+r_2)^{\gamma}}}
k^{\frac{(2r_1-3/2)p}{p-1}-\frac{(2r_2-3/2)p}{\gamma(p-1)}-1/\gamma}\right)^{(p-1)/p}
$$
$$
\leq c \delta n^{2r_1-3/2+\frac{p-1}{p}} \, \mathop{{\sum}}\limits_{l_{r_1}=0}^{n/r_2^\gamma-r_1} (l_{r_1}+1/2)\,
$$
\begin{equation} \label{Sec_dif_C}
\times \mathop{{\sum}}\limits_{m_{r_2}=0}^{(\frac{n}{l_{r_1}+r_1})^{1/\gamma}-r_2}
(m_{r_2}+1/2)^{1-\gamma(2r_1-3/2)+2r_2-3/2+\frac{p-1}{p}-\frac{\gamma(p-1)}{p}} .
\end{equation}
 Whence it follows for $r_1= r_2+1$ and $\gamma= 1$ that
$$
\|\mathcal{D}_{n,1}^{(r_1,r_1-1)} f - \mathcal{D}_{n,1}^{(r_1,r_1-1)} f_\delta\|_{C} \leq c \delta
n^{2r_1-3/2+\frac{p-1}{p}} \, \ln n \, \mathop{{\sum}}\limits_{l_{r_1}=0}^{n/r_2^\gamma-r_1} (l_{r_1}+1/2) = c \delta
n^{2r_1-1/p+3/2} \ln n ,
$$
which was required to prove.

Let's now $r_1 = r_2+1$, $\gamma>1$ or $r_1 > r_2+1$,  $\gamma\geq 1$. Then from (\ref{Sec_dif_C}) we obtain the
desired estimate
$$
\|\mathcal{D}_{n,\gamma}^{(r_1,r_2)} f - \mathcal{D}_{n,\gamma}^{(r_1,r_2)} f_\delta\|_{C} \leq c\delta
n^{2r_1-3/2+\frac{p-1}{p}} \, \mathop{{\sum}}\limits_{l_{r_1}=0}^{n/r_1-r_1} (l_{r_1}+1/2) \leq  c \delta
n^{2r_1-1/p+3/2} .
$$
To estimate $\|\mathcal{D}_{n,\gamma}^{(r_1,r_2)} f - \mathcal{D}_{n,\gamma}^{(r_1,r_2)} f_\delta\|_{C}$ for $r_1=r_2$,
it is enough to repeat an analogous line of reasoning.

In the case of $p=1$ and $p=\infty$ the assertion of Lemma is proved similarly. \vspace{0.1in}

\vskip -4mm

  ${}$ \ \ \ \ \ \ \ \ \ \ \ \ \ \ \ \ \ \ \ \ \ \ \ \ \ \ \ \ \ \ \ \ \ \ \ \ \ \
  \ \ \ \ \ \ \ \ \ \ \ \ \ \ \ \ \ \ \ \ \ \
  \ \ \ \ \ \ \ \ \ \ \ \ \ \ \ \ \ \ \ \ \ \
  \ \ \ \ \ \ \ \ \ \ \ \ \ \ \ \ \ \ \ \ \ \
  \ \ \ \ \ \ \ \ \ \ \ \ \ \ \ \ \ \ \ \ \ \
  \ \ \ \ \ \ \ \ \ \ \ \ \ \ \ \ \ \ \ \ \ \
  \ \ $\Box$

The combination of Lemmas \ref{lemma_BoundErrHCC} and \ref{lemma_BoundPertHCC} gives
\begin{theorem} \label{Th2}
Let $f\in L^\mu_{s,2}$, $1\leq s< \infty$, $r_1 \geq r_2$, $\mu>2r_1-1/s+3/2$, and let the condition
(\ref{perturbation}) be satisfied.

I)  If $r_1 = r_2$, then for $n\asymp \left(\delta/ \ln^{1/p-1/s} \frac{1}{\delta}\right)^{-\frac{1}{\mu-1/p+1/s}}$ and
$\gamma=1$ it holds true
  $$
   \|f^{(r_1,r_1)} - \mathcal{D}^{(r_1,r_1)}_{n,1} f_\delta\|_{C} \leq c \left(\delta \ln^{1/s-1/p}
\frac{1}{\delta}\right)^{\frac{\mu-2r_1+1/s-3/2}{\mu-1/p+1/s}} \ln^{2-1/s} \frac{1}{\delta} .
  $$

II)  If $r_1 = r_2+1$, then\\
a) for $n\asymp \delta^{-\frac{1}{\mu-1/p+1/s}}$ and $\gamma=1$ it holds true
$$
  \|f^{(r_1,r_1-1)}-\mathcal{D}_{n,1}^{(r_1,r_1-1)} f_\delta\|_{C}\leq c \delta^{\frac{\mu-2r_1+1/s-3/2}{\mu-1/p+1/s}} \ln \frac{1}{\delta}  ,
$$
b) for $n\asymp \delta^{-\frac{1}{\mu-1/p+1/s}}$ and any \,$ \gamma \in \left(1,
\frac{\mu-2r_1+1/s+5/2}{\mu-2r_1+1/s+1/2}\right) \bigcup \left( \frac{\mu-2r_1+1/s+5/2}{\mu-2r_1+1/s+1/2},
\frac{\mu-2r_1+1/s+1/2}{\mu-2r_1+1/s-3/2}\right) $ \\ it holds true
$$
  \|f^{(r_1,r_1-1)}-\mathcal{D}_{n,\gamma}^{(r_1,r_1-1)} f_\delta\|_{C}\leq c \delta
^{\frac{\mu-2r_1+1/s-3/2}{\mu-1/p+1/s}} ,
$$
c) for $n\asymp \left(\delta/ \ln \frac{1}{\delta}\right)^{-\frac{1}{\mu-1/p+1/s}}$ and $\, \gamma =
\frac{\mu-2r_1+1/s+1/2}{\mu-2r_1+1/s-3/2}\, $ it holds true
$$
  \|f^{(r_1,r_1-1)}-\mathcal{D}_{n,\gamma}^{(r_1,r_1-1)} f_\delta\|_{C}\leq c \left(\delta/ \ln
\frac{1}{\delta}\right)^{\frac{\mu-2r_1+1/s-3/2}{\mu-1/p+1/s}} \ln \frac{1}{\delta}  ,
$$
d) for $n\asymp \left(\delta/ \ln^{1-1/s} \frac{1}{\delta}\right)^{-\frac{1}{\mu-1/p+1/s}}$ and $\, \gamma =
\frac{\mu-2r_1+1/s+5/2}{\mu-2r_1+1/s+1/2}\, $ it holds true
$$
  \|f^{(r_1,r_1-1)}-\mathcal{D}_{n,\gamma}^{(r_1,r_1-1)} f_\delta\|_{C}\leq c \left(\delta/ \ln^{1-1/s}
\frac{1}{\delta}\right)^{\frac{\mu-2r_1+1/s-3/2}{\mu-1/p+1/s}} \ln^{1-1/s} \frac{1}{\delta}  .
$$

III)  If $r_1 > r_2+1$, then\\
a) for $n\asymp \delta^{-\frac{1}{\mu-1/p+1/s}}$ and any \, $ \gamma \in \left[1,
\frac{\mu-2r_2+1/s-3/2}{\mu-2r_1+1/s+1/2}\right) \bigcup \left(
\frac{\mu-2r_2+1/s-3/2}{\mu-2r_1+1/s+1/2}, \frac{\mu-2r_2+1/s+1/2}{\mu-2r_1+1/s+1/2}\right) $\\
$
 \bigcup \left(
\frac{\mu-2r_2+1/s+1/2}{\mu-2r_1+1/s+1/2}, \frac{\mu-2r_2+1/s-3/2}{\mu-2r_1+1/s-3/2}\right) $ \, it holds true
$$
  \|f^{(r_1,r_2)}-\mathcal{D}_{n,\gamma}^{(r_1,r_2)} f_\delta\|_{C}\leq c \delta
^{\frac{\mu-2r_1+1/s-3/2}{\mu-1/p+1/s}} ,
$$
b) for $n\asymp \left(\delta/ \ln  \frac{1}{\delta}\right)^{-\frac{1}{\mu-1/p+1/s}}$ and $\, \gamma =
\frac{\mu-2r_2+1/s-3/2}{\mu-2r_1+1/s+1/2}, \quad \gamma = \frac{\mu-2r_2+1/s-3/2}{\mu-2r_1+1/s-3/2}\, $ it holds true
$$
  \|f^{(r_1,r_2)}-\mathcal{D}_{n,\gamma}^{(r_1,r_2)} f_\delta\|_{C}\leq c \left(\delta/ \ln
\frac{1}{\delta}\right)^{\frac{\mu-2r_1+1/s-3/2}{\mu-1/p+1/s}} \ln \frac{1}{\delta}  ,
$$
c) for $n\asymp \left(\delta/ \ln^{1-1/s} \frac{1}{\delta}\right)^{-\frac{1}{\mu-1/p+1/s}}$ and $\, \gamma =
\frac{\mu-2r_2+1/s+1/2}{\mu-2r_1+1/s+1/2}\, $ it holds true
$$
  \|f^{(r_1,r_2)}-\mathcal{D}_{n,\gamma}^{(r_1,r_2)} f_\delta\|_{C}\leq c \left(\delta/ \ln^{1-1/s}
\frac{1}{\delta}\right)^{\frac{\mu-2r_1+1/s-3/2}{\mu-1/p+1/s}} \ln^{1-1/s} \frac{1}{\delta}  .
$$
\end{theorem}

\vskip 2mm

\begin{remark} \label{Rem3}
 \rm The $\mathcal{D}^{(r,r)}_{n}$ (\ref{ModVer}) method was studied earlier (see \cite{Sem_Sol_2021}) for the problem of
 numerical differentiation of functions from
 $L^{\mu}_{s,2}$  in the case of $r=1$ and $p=s=2$. Thus, the results of theorems \ref{Th1} and \ref{Th2} generalize
 studies \cite{Sem_Sol_2021} for the case of arbitrary $r,p,s$.
\end{remark}

\vskip 2mm
\section{Minimal radius of Galerkin information}

Let us turn to finding sharp estimates (in the power scale) for the minimal radius. First, we establish a lower
estimate for the quantity $R_{N,\delta}^{(r_1,r_2)}(L^{\mu}_{s,2}, C, \ell_p)$. We fix an arbitrarily chosen domain
 $\hat{\Omega}$, $\card(\hat{\Omega})\leq N$, of the coordinate plane $[r_1,\infty)\times[r_2,\infty)$ and build an auxiliary function
 $$
 f_1(t,\tau)  = \widetilde{c} \, \bigg( \varphi_0(t) \varphi_0(\tau) \,
 + \,  N^{-\mu-1/s}\, r_2^{-\mu} \varphi_{r_2}(\tau) \mathop{{\sum}\, '}\limits_{k=N+r_1}^{3N+r_1}
 \varphi_k(t) \bigg) ,
 $$
 where the sum $\mathop{{\sum}\, '}\limits_{k=N+r_1}^{3N+r_1}$ is taken over any $N$ pairwise distinct functions $\varphi_k(t)$ such that
 $N+r_1\leq k\leq 3N+r_1$ and $(k,r_2)\notin \hat{\Omega}$. Obviously, there  is at least one set of such functions.

Now we estimate the norm of $f_1$ in the space metric  $L^{\mu}_{s,2}$:
 $$
 \|f_1\|^s_{s,\mu} = {\widetilde{c}}^{\,s} \, \bigg( 1
 +  N^{-s\mu-1}  \,  \mathop{{\sum}\, '}\limits_{k=N+r_1}^{3N+r_1}
 k^{s\mu} \bigg) \leq
 {\widetilde{c}}^{\,s}    \,  \bigg( 1
 + 4^{s\mu} \bigg) .
 $$
 Whence it follows that to satisfy the condition $\|f_1\|_{s,\mu}\leq 1$ it suffices to take
 $$
 {\widetilde{c}} =  \bigg( 1
 + 4^{s\mu} \bigg)^{-1/s} .
 $$

Next, we take another function from the class $L^{\mu}_{s,2}$:
 $$
 f_2(t,\tau) = \widetilde{c} \, \varphi_0(t) \varphi_0(\tau)  .
 $$

Let us find a lower bound for the quantity $\|f_1^{(r_1,r_2)}-f_2^{(r_1,r_2)}\|_{C}$. For this we need formulas
$$
\varphi_r^{(r)}(t) = \frac{\sqrt{r+1/2}}{2^{r-1/2}}\ \frac{(2r)!}{r!}\ \varphi_0(t), \quad f_2^{(r_1,r_2)}(t,\tau)
\equiv 0 ,
$$
 $$
 f_1^{(r_1,r_2)}(t,\tau) = \frac{\widetilde{c}}{r_2^\mu} \, {N^{-\mu-1/s}} \,
\varphi_{r_2}^{(r_2)}(\tau) \mathop{{\sum}\, '}\limits_{k=N+r_1}^{3N+r_1} \varphi_k^{(r_1)}(t)
$$
$$
= 2^{r_1}\, \frac{\widetilde{c}}{r_2^\mu} \, {N^{-\mu-1/s}} \, \varphi_{r_2}^{(r_2)}(\tau) \mathop{{\sum}\,
'}\limits_{k=N+r_1}^{3N+r_1} \sqrt{k+1/2} \mathop{{\sum}^*}\limits_{l_1=r_1-1}^{k-1} (l_1+1/2)
\mathop{{\sum}^*}\limits_{l_2=r_1-2}^{l_1-1} (l_2+1/2)
$$
\begin{equation}  \label{f_1^(r,r)}
\ldots \mathop{{\sum}^*}\limits_{l_{r_1-1}=1}^{l_{r_1-2}-1} (l_{r_1-1}+1/2)
\mathop{{\sum}^*}\limits_{l_{r_1}=0}^{l_{r_1-1}-1} \sqrt{l_{r_1}+1/2}\, \varphi_{l_{r_1}}(t) .
\end{equation}
It is easy to see that
$$
 \|f_1^{(r_1,r_2)}-f_2^{(r_1,r_2)}\|_{C}
 \geq  |f_1^{(r_1,r_2)}(1,1)|
 \geq  \overline{c} \, N^{-\mu+2r_1-1/s+3/2}  ,
$$
where
 $$
\overline{c} = \frac{\widetilde{c}\sqrt{r_2+1/2}}{2^{r_1+r_2} r_2^\mu}\,  \frac{(2r_2)!}{r_1! r_2!} .
 $$
Since for any $1\leq p\leq \infty$ it holds true
 $$
 \|\overline{f}_1-\overline{f}_2\|_{\ell_p}
 = \frac{\widetilde{c}}{r_2^\mu} \, N^{-\mu-1/s+1/p} ,
 $$
 then in the case of $N^{-\mu-1/s+1/p}\leq r_2^\mu \delta/{\widetilde{c}}$  under $\delta$-perturbations of the functions $f_1$ and $f_2$
 can be considered
$$
f^\delta_1(t,\tau) = f_2(t,\tau), \qquad f^\delta_2(t,\tau) = f_1(t,\tau) .
$$

Let us find  the upper bound  for $\|f_1^{(r_1,r_2)}-f_2^{(r_1,r_2)}\|_{C}$.
 Taking into account the relation $G(\hat{\Omega},\overline{f}_1^{\delta})=G(\hat{\Omega},\overline{f}_2^{\delta})$,
for any $\psi^{(r_1,r_2)}(\hat{\Omega})\in\Psi(\hat{\Omega})$ we find
 $$
 \|f_1^{(r_1,r_2)}-f_2^{(r_1,r_2)}\|_{C}
 \leq \|f_1^{(r_1,r_2)}-\psi^{(r_1,r_2)}(G(\hat{\Omega},\overline{f}_1^{\delta}))\|_{C}
 + \|f_2^{(r_1,r_2)}-\psi^{(r_1,r_2)}(G(\hat{\Omega},\overline{f}_2^{\delta}))\|_{C}
 \leq
 $$
 $$
 \leq 2 \, \sup\limits_{f\in L^{\mu}_{s,2}, \|f\|_{s,\mu}\leq 1}
 \ \sup\limits_{\overline{f^\delta}: \, (\ref{perturbation})}
 \| f^{(r_1,r_2)} - \psi^{(r_1,r_2)}(G(\hat{\Omega},\overline{f}^{\delta})) \|_C
 =: 2 \, \varepsilon_{\delta}(L^{\mu}_{s,2}, \psi^{(r_1,r_2)}(\hat{\Omega}), C, \ell_p) .
 $$

That is
 $$
 \varepsilon_{\delta}(L^{\mu}_{s,2}, \psi^{(r_1,r_2)}(\hat{\Omega}), C, \ell_p)
 \geq \frac{\overline{c}}{2} \, N^{-\mu+2r_1-1/s+3/2} .
 $$

From the fact that the domain $\hat{\Omega}$
 and the algorithm $\psi^{(r_1,r_2)}(\hat{\Omega})\in\Psi(\hat{\Omega})$  are arbitrary, follows that
 $$
R_{N,\delta}^{(r_1,r_2)}(L^{\mu}_{s,2}, C, \ell_p)
 \geq \frac{\overline{c}}{2} \, N^{-\mu+2r_1-1/s+3/2} .
 $$

Thus, the following assertion is proved.

\begin{theorem} \label{Th5.1}
Let $1\leq s< \infty$, $r_1\geq r_2$, $\mu>2r_1-1/s+3/2$, $1\leq p \leq \infty$. Then for any
 $N\geq\Big(r_2^\mu\delta/\widetilde{c}\Big)^{-1/(\mu+1/s-1/p)}$ it holds true
 $$
R_{N,\delta}^{(r_1,r_2)}(L^{\mu}_{s,2}, C, \ell_p)
 \geq \frac{\overline{c}}{2} \, N^{-\mu+2r_1-1/s+3/2} .
 $$
\end{theorem}

The following assertion contains sharp estimates (in the power scale) for the minimal radius in the uniform metric.

\begin{theorem} \label{Th5.2}
Let $f\in L^\mu_{s,2}$, $1\leq s< \infty$, $r_1 \geq r_2$, $\mu>2r_1+1/2-3/s$.

I)  If $r_1 = r_2$, then for $N\asymp \left(\delta/ \ln^{\mu} \frac{1}{\delta}\right)^{-\frac{1}{\mu-1/p+1/s}}$ it
holds true
$$
 N^{-\mu+2r_1-1/s+3/2} \preceq
R_{N,\delta}^{(r_1,r_1)}(L^{\mu}_{s,2}, C, \ell_p)
 \preceq N^{-\mu+2r_1-1/s+3/2} \ln^{\mu-2r_1+1/2} N ,
$$
$$
\left(\delta/ \ln^{\mu} \frac{1}{\delta}\right)^{\frac{\mu-2r_1+1/s-3/2}{\mu-1/p+1/s}} \preceq
R_{N,\delta}^{(r_1,r_1)}(L^{\mu}_{s,2}, C, \ell_p) \preceq \left(\delta \ln^{1/s-1/p}
\frac{1}{\delta}\right)^{\frac{\mu-2r_1+1/s-3/2}{\mu-1/p+1/s}} \ln^{2-1/s} \frac{1}{\delta} .
$$
The upper bounds are implemented by the method $\mathcal{D}^{(r_1,r_1)}_{n,\gamma}$ (\ref{ModVer}) with  $n\asymp
\left(\delta/ \ln^{1/p-1/s} \frac{1}{\delta}\right)^{-\frac{1}{\mu-1/p+1/s}}$ and $\gamma=1$.

II)  If $r_1 = r_2+1$, then
for $N\asymp \delta^{-\frac{1}{\mu-1/p+1/s}}$ it holds true
$$
R_{N,\delta}^{(r_1,r_1-1)}(L^{\mu}_{s,2}, C, \ell_p)
 \asymp N^{-\mu+2r_1-1/s+3/2} \asymp \delta^{\frac{\mu-2r_1+1/s-3/2}{\mu-1/p+1/s}} .
$$
The order-optimal bounds are implemented by the method $\mathcal{D}^{(r_1,r_1-1)}_{n,\gamma}$ (\ref{ModVer}) with
$n\asymp \delta^{-\frac{1}{\mu-1/p+1/s}}$ and any\, $ \gamma \in \left(1,
\frac{\mu-2r_1+1/s+5/2}{\mu-2r_1+1/s+1/2}\right) \bigcup \left( \frac{\mu-2r_1+1/s+5/2}{\mu-2r_1+1/s+1/2},
\frac{\mu-2r_1+1/s+1/2}{\mu-2r_1+1/s-3/2}\right) $ .

III)  If $r_1 > r_2+1$, then
for $N\asymp \delta^{-\frac{1}{\mu-1/p+1/s}}$ it holds true
$$
R_{N,\delta}^{(r_1,r_2)}(L^{\mu}_{s,2}, C, \ell_p)
 \asymp N^{-\mu+2r_1-1/s+3/2} \asymp \delta^{\frac{\mu-2r_1+1/s-3/2}{\mu-1/p+1/s}} .
$$
The order-optimal bounds are implemented by the method $\mathcal{D}^{(r_1,r_2)}_{n,\gamma}$ (\ref{ModVer}) for $n\asymp
\delta^{-\frac{1}{\mu-1/p+1/s}}$ and any\, $ \gamma \in \left(1, \frac{\mu-2r_2+1/s-3/2}{\mu-2r_1+1/s+1/2}\right)
\bigcup \left( \frac{\mu-2r_2+1/s-3/2}{\mu-2r_1+1/s+1/2}, \frac{\mu-2r_2+1/s+1/2}{\mu-2r_1+1/s+1/2}\right)
 \bigcup \left(
\frac{\mu-2r_2+1/s+1/2}{\mu-2r_1+1/s+1/2}, \frac{\mu-2r_2+1/s-3/2}{\mu-2r_1+1/s-3/2}\right) $ .
\end{theorem}

 \bf Proof.
 \rm
The upper bounds for  $R_{N,\delta}^{(r,r)}(L^{\mu}_{s,2}, C, \ell_p)$ follow from Theorem \ref{Th2}.
The lower bound is found in Theorem  \ref{Th5.1}.

\vskip -4mm

  ${}$ \ \ \ \ \ \ \ \ \ \ \ \ \ \ \ \ \ \ \ \ \ \ \ \ \ \ \ \ \ \ \ \ \ \ \ \ \ \
  \ \ \ \ \ \ \ \ \ \ \ \ \ \ \ \ \ \ \ \ \ \
  \ \ \ \ \ \ \ \ \ \ \ \ \ \ \ \ \ \ \ \ \ \
  \ \ \ \ \ \ \ \ \ \ \ \ \ \ \ \ \ \ \ \ \ \
  \ \ \ \ \ \ \ \ \ \ \ \ \ \ \ \ \ \ \ \ \ \
  \ \ \ \ \ \ \ \ \ \ \ \ \ \ \ \ \ \ \ \ \ \
  \ \ $\Box$

 \rm
Let's turn  to estimating the minimal radius in the integral metric.

\begin{theorem} \label{Th5.4}
Let $1\leq s< \infty$, $r_1\geq r_2$, $\mu>2r_1-1/s+1/2$, $1\leq p \leq \infty$. Then for any
$N\geq\Big(r_2^\mu\delta/\widetilde{c}\Big)^{-1/(\mu+1/s-1/p)}$ it holds true
$$
R_{N,\delta}^{(r_1,r_2)}(L^{\mu}_{s,2}, L_2, \ell_p) \geq \overline{\overline{c}} \, N^{-\mu+2r_1-1/s+1/2} ,
$$
 where $\overline{\overline{c}} = \frac{\widetilde{c} \sqrt{r_2+1/2}}{2^{3r_1+r_2-3/2}\, r_2^{\mu}} \frac{(2r_2)!}{(r_1-1)!r_2!}$ .
\end{theorem}

 \rm

 \bf Proof \rm The proof of Theorem \ref{Th5.4} almost completely coincides with the proof of Theorem \ref{Th5.1},
 including the form of the auxiliary functions
 $f_1$, $f_1^{\delta}$, $f_2$, $f_2^{\delta}$.
 The only difference is in the lower estimate of the norm of
 $f_1^{(r_1,r_2)}-f_2^{(r_1,r_2)}$. Changing the order of summation in  (\ref{f_1^(r,r)}) yields to the representation
$$
 f_1^{(r_1,r_2)}(t,\tau) = 2^{r_1}\, \frac{\widetilde{c}}{r_2^\mu} \, {N^{-\mu-1/s}} \, \varphi_{r_2}^{(r_2)}(\tau)
$$
$$
\times \Big(\mathop{{\sum}^*}\limits_{l_{r_1}=0}^{N} \sqrt{l_{r_1}+1/2}\, \varphi_{l_{r_1}}(t) \mathop{{\sum}\,
'}\limits_{k=N+r_1}^{3N+r_1} \sqrt{k+1/2} + \mathop{{\sum}^*}\limits_{l_{r_1}=N+1}^{3N} \sqrt{l_{r_1}+1/2}\,
\varphi_{l_{r_1}}(t) \mathop{{\sum}\, '}\limits_{k=l_{r_1}+r_1}^{3N+r_1} \sqrt{k+1/2}\Big)
$$
$$
\times \mathop{{\sum}^*}\limits_{l_1=l_{r_1}+r_1-1}^{k-1} (l_1+1/2)
\mathop{{\sum}^*}\limits_{l_2=l_{r_1}+r_1-2}^{l_1-1} (l_2+1/2) \ldots
\mathop{{\sum}^*}\limits_{l_{r_1-1}=l_{r_1}+1}^{l_{r_1-2}-1} (l_{r_1-1}+1/2)  .
$$
Let us introduce the following notation
$$
c' = \frac{2^{2(r_1-r_2)+1}}{r_2^{2\mu}}\, \widetilde{c}^2 (r_2+1/2) \frac{((2r_2)!)^2}{(r_2!)^2} .
$$
Then it follows that
 $$
 \|f_1^{(r_1,r_2)}-f_2^{(r_1,r_2)}\|_{L_2}^2 = \|f_1^{(r_1,r_2)}\|_{L_2}^2
 \geq c' \, {N^{-2\mu-2/s}}\,  \mathop{{\sum}^*}\limits_{l_{r_1}=0}^{N} (l_{r_1}+1/2)
 $$
 $$
 \times \left( \mathop{{\sum}\,
'}\limits_{k=N+r_1}^{3N+r_1} \sqrt{k+1/2} \mathop{{\sum}^*}\limits_{l_1=l_{r_1}+r_1-1}^{k-1} (l_1+1/2)
\mathop{{\sum}^*}\limits_{l_2=l_{r_1}+r_1-2}^{l_1-1} (l_2+1/2) \ldots
\mathop{{\sum}^*}\limits_{l_{r_1-1}=l_{r_1}+1}^{l_{r_1-2}-1} (l_{r_1-1}+1/2)\right)^2
$$
$$
  \geq c' \, \frac{N^{-2\mu-2/s} \cdot N^{4r_1-1}}{4^{4(r_1-1)}((r_1-1)!)^2}
 \mathop{{\sum}^*}\limits_{l_{r_1}=0}^{N/2} (l_{r_1}+1/2) .
 $$
 That is, the estimate holds
 $$
 \|f_1^{(r_1,r_2)}-f_2^{(r_1,r_2)}\|_{L_2}
 \geq \frac{\sqrt{c'}}{4^{2r_1-1}(r_1-1)!} N^{-\mu+2r_1-1/s+1/2} .
 $$
 Whence we obtain that the relation
 $$
 \varepsilon_{\delta}(L^{\mu}_{s,2}, \psi^{(r_1,r_2)}(\hat{\Omega}), L_2, \ell_p)
 \geq \overline{\overline{c}} \, N^{-\mu+2r_1-1/s+1/2}
 $$
is true for any $N\geq\Big(r_2^\mu\delta/\widetilde{c}\Big)^{-1/(\mu+1/s-1/p)}$. From the fact that the domain
$\hat{\Omega}$ and the algorithm $\psi^{(r_1,r_2)}(\hat{\Omega})\in\Psi(\hat{\Omega})$  are arbitrary, it follows that
 $$
R_{N,\delta}^{(r_1,r_2)}(L^{\mu}_{s,2}, L_2, \ell_p)
 \geq \overline{\overline{c}} \, N^{-\mu+2r_1-1/s+1/2} .
 $$
 Thus, the proof of Theorem \ref{Th5.4} has been complied.

\vskip -4mm

  ${}$ \ \ \ \ \ \ \ \ \ \ \ \ \ \ \ \ \ \ \ \ \ \ \ \ \ \ \ \ \ \ \ \ \ \ \ \ \ \
  \ \ \ \ \ \ \ \ \ \ \ \ \ \ \ \ \ \ \ \ \ \
  \ \ \ \ \ \ \ \ \ \ \ \ \ \ \ \ \ \ \ \ \ \
  \ \ \ \ \ \ \ \ \ \ \ \ \ \ \ \ \ \ \ \ \ \
  \ \ \ \ \ \ \ \ \ \ \ \ \ \ \ \ \ \ \ \ \ \
  \ \ \ \ \ \ \ \ \ \ \ \ \ \ \ \ \ \ \ \ \ \
  \ \ $\Box$

The following assertion contains sharp estimates (in the power scale) for the minimal radius in the integral metric.
 \rm

\begin{theorem} \label{Th5.5}
Let $f\in L^\mu_{s,2}$, $1\leq s< \infty$, $r_1 \geq r_2$, $\mu>2r_1+1/2-1/s$.

I)  If $r_1 = r_2$, then for $N\asymp \left(\delta/ \ln^{\mu} \frac{1}{\delta}\right)^{-\frac{1}{\mu-1/p+1/s}}$ it
holds true
$$
 N^{-\mu+2r_1-1/s+1/2} \preceq
R_{N,\delta}^{(r_1,r_1)}(L^{\mu}_{s,2}, L_2, \ell_p)
 \preceq N^{-\mu+2r_1-1/s+1/2} \ln^{\mu-2r_1+1} N ,
$$
$$
\left(\delta/ \ln^{\mu} \frac{1}{\delta}\right)^{\frac{\mu-2r_1+1/s-1/2}{\mu-1/p+1/s}} \preceq
R_{N,\delta}^{(r_1,r_1)}(L^{\mu}_{s,2}, L_2, \ell_p) \preceq \left(\delta \ln^{1/s-1/p}
\frac{1}{\delta}\right)^{\frac{\mu-2r_1+1/s-1/2}{\mu-1/p+1/s}} \ln^{3/2-1/s} \frac{1}{\delta} .
$$
The upper bounds are implemented by the method $\mathcal{D}^{(r_1,r_1)}_{n,\gamma}$ (\ref{ModVer}) with  $n\asymp
\left(\delta^{-1} \ln^{1/p-1/s} \frac{1}{\delta}\right)^{\frac{1}{\mu-1/p+1/s}}$ and $\gamma=1$.

II)  If $r_1 > r_2$, then
for $N\asymp \delta^{-\frac{1}{\mu-1/p+1/s}}$ it holds true
$$
R_{N,\delta}^{(r_1,r_2)}(L^{\mu}_{s,2}, L_2, \ell_p)
 \asymp N^{-\mu+2r_1-1/s+1/2} \asymp \delta^{\frac{\mu-2r_1+1/s-1/2}{\mu-1/p+1/s}} .
$$
The order-optimal bounds are implemented by the method $\mathcal{D}^{(r_1,r_2)}_{n,\gamma}$ (\ref{ModVer}) with
$n\asymp \delta^{-\frac{1}{\mu-1/p+1/s}}$ and any\, $ \gamma \in \left(1,
\frac{\mu-2r_2+1/s-1/2}{\mu-2r_1+1/s+1/2}\right) \bigcup \left( \frac{\mu-2r_2+1/s-1/2}{\mu-2r_1+1/s+1/2},
\frac{\mu-2r_2+1/s+1/2}{\mu-2r_1+1/s+1/2}\right)
 \bigcup \left(
\frac{\mu-2r_2+1/s+1/2}{\mu-2r_1+1/s+1/2}, \frac{\mu-2r_2+1/s-1/2}{\mu-2r_1+1/s-1/2}\right) $ .
\end{theorem}

{\textbf{\textit{Proof.}}
 \rm The upper bounds for  $R_{N,\delta}^{(r_1,r_2)}(L^{\mu}_{s,2}, L_2, \ell_p)$ follow from Theorem
\ref{Th1}.
The lower bound is found in Theorem \ref{Th5.4}.

\vskip -4mm

  ${}$ \ \ \ \ \ \ \ \ \ \ \ \ \ \ \ \ \ \ \ \ \ \ \ \ \ \ \ \ \ \ \ \ \ \ \ \ \ \
  \ \ \ \ \ \ \ \ \ \ \ \ \ \ \ \ \ \ \ \ \ \
  \ \ \ \ \ \ \ \ \ \ \ \ \ \ \ \ \ \ \ \ \ \
  \ \ \ \ \ \ \ \ \ \ \ \ \ \ \ \ \ \ \ \ \ \
  \ \ \ \ \ \ \ \ \ \ \ \ \ \ \ \ \ \ \ \ \ \
  \ \ \ \ \ \ \ \ \ \ \ \ \ \ \ \ \ \ \ \ \ \
  \ \ $\Box$

 \section{ Acknowledgements }
This project has received funding through the MSCA4Ukraine project,
which is funded by the European Union.
In addition,
the first named author is supported by the Volkswagen Foundation project "From Modeling and Analysis to Approximation".
Also, the authors acknowledge partial financial support due to the project "Mathematical modelling of complex dynamical systems and processes caused by the state security" (Reg. No. 0123U100853).

\begin{small}
\begin{flushleft}
\textsc{Institute of Mathematics, National Academy of Sciences of Ukraine,
3, Tereschenkivska Str., 01024, Kiev, Ukraine}\\
\textit{E-mail address}: semenovaevgen@gmail.com, solodky@imath.kiev.ua
\end{flushleft}
\end{small}
\begin{flushright}
Received 11.02.2009; revised 23.03.2009
\end{flushright}


\begin{thebibliography}{99}
    \begin{small}

        \bibitem{Ahn&Choi&Ramm_2006}
        S.~Ahn, U.J.~Choi,  A.G.~Ramm, A scheme for stable numerical differentiation
        // J. Comput. Appl. Math. -- 2006.-- Vol.~
        186(2).-- P. 325--334


        \bibitem{And84} R. S. Anderssen, F. R. de Hoog, Finite difference
        methods for the numerical differentiation of non-
        exact data, Computing, 33(1984), 259–267.


        \bibitem{Cul71} Cullum, J., Numerical Differentiation and Regularization. SIAM Journal on Numerical Analysis.--
        1971. -- 8(2), 254–265. http://www.jstor.org/stable/2949474

        \bibitem{Dolgopolova&Ivanov_USSR_Comput_Math_Math_Phys_1966_Eng}
        T.F.~Dolgopolova, V.K.~Ivanov,
        On numerical differentiation
        // Zh. Vychisl. Mat.and Mat. Ph.-- 1966. -- Vol.~6(3).-- P. 223--232.

        \bibitem{EgorKond_1989}
        Yu.V.Egorov, V.A.Kondrat'ev, On a problem of numerical differentiation // Vestnik Moskov. Univ. Ser. I Mat. Mekh. 3 (1989), 80-81.

        \bibitem{ErbSem2015} Wolfgang Erb, \, Evgeniya V. Semenova,
        On adaptive discretization schemes for the solution of ill-posed problems with semiiterative methods// Applicable
        Analysis, Vol.94, 2015(10)

        \bibitem{Groetsch_1992_V74_N2}
        Groetsch~C.W.,
        Optimal order of accuracy in Vasin's method for differentiation of noisy functions
        //J. Optim.Theory Appl.-- 1992.-- Vol.~74(2).-- P. 373--378.

        \bibitem{Hanke&Scherzer_2001_V108_N6}
        M.~Hanke, O.~Scherzer,
        Inverse problems light: numerical differentiation
        //  Amer. Math. Monthly.-- 2001.-- Vol.~108(6).-- P. 512--521.

        \bibitem{Lu&Naum&Per}
        S. Lu, V. Naumova, S.\, V. Pereverzev,
        Legendre polynomials as a recommended basis for numerical differentiation in the
        presence of stochastic white noise
        // J. Inverse Ill-Posed Probl.-- 2013.-- Vol.~21(2).-- P. 193--216.

        \bibitem{Meng&Zhaoa&Mei&Zhou_2020}
        Z.~Meng, Z.~Zhaoa, D.~Mei, Y.~Zhou,
        Numerical differentiation for two-dimensional functions by a Fourier extension method
        // Inverse Problems in Science and Engineering.-- 2020.-- Vol.~28(1).-- P. 1--18.

        \bibitem{Wang_Hon_Ch_2006}Y.B. Wang, Y.C. Hon, J. Cheng
        Reconstruction of high order derivatives from input data
        J. Inverse Ill-posed Probl., 14 (2006), pp. 205-218

        \bibitem{Nakamura&Wang&Wang_2008}
        G. Nakamura, S.\,Z. Wang, Y.\,B. Wang, Numerical differentiation for the second order derivatives of functions of two
        variables
        // J. Comput. Appl. Math.-- 2008.-- Vol.  212(2).-- P. 341--358.


        \bibitem{Mileiko_Solodkii_2014}
        G.L.~Myleiko, S.G.~Solodky, The minimal radius of Galerkin information for severely ill-posed problems
        //  Journal of
        Inverse and Ill-Posed Problems.-- 2014.-- Vol.~22(5).-- P. 739--757.

        \bibitem{Mileiko_Solodkii_2016_ApAn}
        G.L.~Mileyko, S.G.~Solodky, On optimization of projection methods for solving some classes of severely ill-posed
        problems // Applicable Analysis. -- 2016. Vol.95(4), P.826-841.

        \bibitem{Mileiko_Solodkii_2017_UMJ}
        G.L.~Mileyko, S.G.~Solodky,
        Hyperbolic cross and complexity of different classes of linear ill-posed problems
        // Ukr. Mat. J.-- 2017.-- Vol.~69(7).-- P. 951--963.

        \bibitem{Mul69}
        M\" uller~C., \
        Foundations of the Mathematical Theory of Electromagnetic Waves
        /C.~ M\" uller --
        Springer--Verlag, Berlin,
        Heidelberg, New York, 1969.


        \bibitem{Pereverzev_Computing_1995}  S.\,V. Pereverzev,
        Optimization of projection methods for solving ill-posed problems, // Computing. -- 1995 -- Vol.55(2), P.113--124.

        \bibitem{PS1996}
        S.\,V. Pereverzev, S.\,G. Solodky, The minimal radius of Galerkin information for the Fredholm problem of the first
        kind, {\it Journal of Complexity} {\bf 12}~(4) (1996), 401--415.

        \bibitem{Qian&Fu&Xiong&Wei_2006}
        Qian~Z.,\
        Fourier truncation method for high order numerical derivatives
        /Z.~Qian,  C.L.~Fu,  X.T.~Xiong,  T.~Wei,
        // Appl.
        Math. Comput.-- 2006.-- Vol.~181(2).-- P. 940--948.

        \bibitem{Qu96} Qu, Ruibin,
        A new approach to numerical differentiation and integration. Mathematical and Computer Modelling 24 (1996): 55-68.

        \bibitem{Ramm_1968_No11}
        Ramm ~ A.G.,
        On numerical differentiation
        // Izv. Vuzov. Matem.-- 1968.-- Vol.~11.-- P. 131--134.

        \bibitem{RammSmir_2001}
        A.G. Ramm, A.B. Smirnova,
        On stable numerical differentiation Math. Comput., 70 (2001), pp. 1131-1153

        \bibitem{Sem_Sol_2021}
        Semenova Y.V., Solodky S.G., Error bounds for Fourier-Legendre truncation method in numerical differentiation, /
        // Journal of Numerical and Applied Mathematics.-- 2021.-- Vol.~137(3).--P.113--130.

        \bibitem{SSS_CMAM}
        Semenova, Evgeniya V., Solodky, Sergiy G. and Stasyuk, Serhii A., Application of Fourier Truncation Method to Numerical
        Differentiation for Bivariate Functions// Computational Methods in Applied Mathematics, vol. 22, no. 2, 2022, pp.
        477-491. https://doi.org/10.1515/cmam-2020-0138

        \bibitem{SSS_Rew2021} Semenova Y., Solodky S.G., Stasyuk S., Truncation method for numerical diffirentiation problem.
        Proceedings of the Institute of Mathematics of the National Academy of Sciences of Ukraine, Vol. 18 No. 1 (2021):
        Modern problems of mathematics and its applications, II

        \bibitem{Sol_Stas_JC2020} S.\,G. Solodky, S.\,A. Stasyuk, Estimates of efficiency for two methods of stable numerical summation of
        smooth functions, {\it Journal of Complexity} {\bf 56} (2020) \\  https://doi.org/10.1016/j.jco.2019.101422

        \bibitem{Sol_Stas_UMZ2022}
        Solodky S.G., Stasyuk S.,
        On optimization of methods of numerical differentiation for bivariate functions
        // Ukr. Mat. J.-- 2022.-- Vol.~74(2).-- P. 253--273.

        \bibitem{TrWW} J.\,F. Traub, G.\,W. Wasilkowski, H. Wozniakowski, {\it Information-Based Complexity}, Academic
        Press, New York, 1988.

        \bibitem{TrW} J.\,F. Traub,  H. Wozniakowski, A General Theory of Optimal Algorithms, Academic Press, New York,
        1980.

        \bibitem{VasinVV_1969_V7_N2}
        Vasin~V.V., \
        Regularization of the numerical differentiation problem
        /V.V. ~Vasin
        // Mat. app. Ural un-t.-- 1969.-- Vol.~7(2).-- P. 29--33.

        \bibitem{Zhao_2010}
        Zhao~Z.,
        A truncated Legendre spectral method for solving numerical differentiation
        //  International Jounal of
        Computer Mathematics.-- 2010.-- Vol.~87.-- P. 3209--3217.

        \bibitem{Zhao&Meng&Zhao&You&Xie_2016}
        Z.~Zhao, Z.~Meng, L.~Zhao, L.~You, O.~Xie,
        A stabilized algorithm for multi-dimensional numerical differentiation
        // Journal of Algorithms and Computational Technology.-- 2016.-- Vol.~10(2).-- P. 73--81.



        %

        %

    \end{small}
\end{thebibliography}
\end{document}